\documentclass[psamsfonts]{jams-l}
\usepackage{amsmath,amssymb,amsthm,xspace,url,graphicx,hyperref}

\newtheorem{theorem}{Theorem}[section]
\newtheorem{lemma}[theorem]{Lemma}
\newtheorem{cor}[theorem]{Corollary}

\theoremstyle{definition}

\theoremstyle{remark}
\newtheorem{remark}[theorem]{Remark}

\numberwithin{equation}{section}
 \renewcommand{\MRhref}[2]{\href{http://www.ams.org/mathscinet-getitem?mr=#1}{MR#2}}

\makeatletter
\def\@strippedMR{}
\def\@scanforMR#1#2#3\endscan{%
  \ifx#1M\ifx#2R\def\@strippedMR{#3}%
  \else\def\@strippedMR{#1#2#3}%
  \fi\fi}
\def\@rst #1 #2other{#1}
\renewcommand\MR[1]{\relax\ifhmode\unskip\spacefactor3000 \space\fi
  \@scanforMR#1\endscan
  \MRhref{\expandafter\@rst \@strippedMR other}{\@strippedMR}}
\newcommand\MRs[1]{\relax\ifhmode\unskip\spacefactor3000 \space\fi
  \@scanforMR#1\endscan
  \MRhref{\@strippedMR}{\@strippedMR}}
\makeatother

\def\rcs $#1: #2 ${\expandafter\def\csname rcs#1\endcsname {#2}}
\rcs $Date: 2008/03/31 12:09:54 $
\rcs $Revision: 1.31 $

\newcommand{\old}[1]{}
\newcommand{\thh}{\ensuremath{^{\text{th}}}\xspace}

\newcommand{\st}{\ensuremath{^{\text{st}}}\xspace}
\newcommand{\R}{{\mathbb R}}
\newcommand{\N}{{\mathbb N}}

\newcommand{\E}{{\mathbb E}}

\newcommand{\C}{{\mathbb C}}
\newcommand{\Z}{{\mathbb Z}}

\newcommand{\supp}{{\operatorname{supp}}}
\newcommand{\sign}{{\operatorname{sign}}}

\newcommand{\Lip}{{\operatorname{Lip}}}
\newcommand{\diam}{{\operatorname{diam}}}

\newcommand{\I}{{\rm I}\xspace}
\newcommand{\II}{{\rm II}\xspace}
\newcommand{\favored}{{\mathfrak F}}
\newcommand{\eps}{\varepsilon}

\renewcommand{\Re}{\operatorname{Re}}

\renewcommand{\phi}{\varphi}
\newcommand{\increasing}{\mbox{$\star$-increasing}\xspace}
\newcommand{\decreasing}{\mbox{$\star$-decreasing}\xspace}
\renewcommand{\setminus}{\smallsetminus}
\renewcommand{\vec}{}

\begin{document}
\title{Tug-of-war and the infinity Laplacian}

\author{Yuval Peres}
\address{\vspace*{-7pt}Microsoft Research; One Microsoft Way; Redmond, WA 98052}
\address{\mbox{Department of Statistics; 367 Evans Hall; University of California; Berkeley, CA 94720}}
\thanks{Research of Y.P. and S.S. supported in part by
 NSF grants DMS-0244479 and DMS-0104073.}
\urladdr{\url{http://research.microsoft.com/~peres/}}

\author{Oded Schramm}
\address{Microsoft Research; One Microsoft Way; Redmond, WA 98052}
\urladdr{\url{http://research.microsoft.com/~schramm/}}

\author{Scott Sheffield}
\address{\vspace*{-7pt}Microsoft Research; One Microsoft Way; Redmond, WA 98052}
\address{\mbox{Department of Statistics; 367 Evans Hall; University of California; Berkeley, CA 94720}}
\curraddr{Courant Institute; 251 Mercer Street; New York, NY 10012}
\urladdr{\url{http://www.cims.nyu.edu/~sheff/}}

\author{David B. Wilson}
\address{Microsoft Research; One Microsoft Way; Redmond, WA 98052}
\urladdr{\url{http://dbwilson.com}}

\subjclass[2000]{91A15, 91A24, 35J70, 54E35, 49N70}
\date{July 8, 2006, in revised form March 31, 2008.}
\keywords{infinity Laplacian, absolutely minimal Lipschitz extension, tug-of-war}
\copyrightinfo{2008}{American Mathematical Society}

\maketitle

\section{Introduction and preliminaries}
\subsection{Overview}
We consider a class of zero-sum two-player stochastic games
called \textbf{tug-of-war} and use them to prove that every
bounded real-valued Lipschitz function $F$ on a subset $Y$ of a
length space $X$ admits a unique \textbf{absolutely minimal (AM)}
extension to $X$, i.e., a unique Lipschitz extension $u:X
\rightarrow \R$ for which $\Lip_Uu = \Lip_{\partial U}u$ for
all open $U \subset X\smallsetminus Y$.  We present an example
that shows this is not generally true when $F$ is merely
Lipschitz and positive.
(Recall that a metric space $(X,d)$ is a \textbf{length space} if
for all $x,y\in X$, the distance $d(x,y)$ is the infimum of the
lengths of continuous paths in $X$ that connect $x$ to $y$. Length
spaces are more general than geodesic spaces, where the infima
need to be achieved.)

When $X$ is the closure of a bounded domain $U \subset \R^n$ and
$Y$ is its boundary, a Lipschitz extension $u$ of $F$ is AM if and
only it is \textbf{infinity harmonic} in the interior of $X
\smallsetminus Y$, i.e., it is a \textbf{viscosity solution}
(defined below) to $\Delta_\infty u = 0$, where $\Delta_\infty$ is
the so-called \textbf{infinity Laplacian}

\begin{equation} \label{pdedef}
\Delta_\infty u = |\nabla u|^{-2} \sum_{i,j} u_{x_i} u_{x_ix_j}
u_{x_j}
\end{equation}
(informally, this is the second derivative of $u$ in the direction
of the gradient of $u$).  Aronsson proved this equivalence for smooth $u$ in 1967
and Jensen proved the general statement in 1993 \cite{MR0217665,
  MR1218686}.  Our analysis of tug-of-war also shows that in this
setting $\Delta_\infty u = g$ has a unique viscosity solution
(extending $F$) when $g:\overline U\to \R$ is continuous and $\inf g > 0$ or $\sup g <
0$, but not necessarily when $g$ assumes values of both signs.
We note that in the study of the homogenous equation
$\Delta_\infty u=0$, the normalizing factor $|\nabla u|^{-2}$
in (\ref{pdedef}) is sometimes omitted; however, it is important to include it in the nonhomogenous equation.
Observe that with the normalization, $\Delta_\infty$ coincides with the ordinary Laplacian $\Delta$ in the
one-dimensional case.

Unlike the ordinary Laplacian or the $p$-Laplacian for $p <
\infty$, the infinity Laplacian can be defined on  \textbf{any}
length space with no additional structure (such as a measure or a
canonical Markov semigroup)---that is, we will see that viscosity
solutions to $\Delta_\infty u = g$ are well-defined in this
generality.  We will establish the above stated uniqueness of solutions $u$ to
$\Delta_\infty u = g$ in the setting of length spaces.

Originally, we were motivated not by the infinity Laplacian but by
random turn Hex~\cite{pssw-hex} and its generalizations,
which led us to consider the tug-of-war game.
As we later learned, tug-of-war games have been considered by
Lazarus, Loeb, Propp and Ullman in~\cite{MR1427981} (see also~\cite{MR1685133}).

Tug-of-war on a metric space is very
natural and conceivably applicable (like differential game theory)
to economic and political modeling.

The intuition provided by thinking of strategies for tug-of-war
yields new results even in the classical setting of domains
in $\R^n$. For instance, in Section~\ref{section:porous}
we show that if $u$ is infinity-harmonic in the unit disk
and its boundary values are in $[0,1]$ and supported on a $\delta$-neighborhood
of the ternary Cantor set on the unit circle, then $u(0)<\delta^\beta$ for some
$\beta>0$.

Before precisely stating our main results, we need several definitions.

\subsection{Random turn games and values}

We consider two-player, zero-sum \textbf{random-turn games}, which
are defined by the following parameters: a set $X$ of
\textbf{states} of the game, two directed \textbf{transition
graphs} $E_{\I}, E_{\II}$ with vertex set $X$, a nonempty set $Y \subset
X$ of \textbf{terminal states} (a.k.a.\ \textbf{absorbing states}),
a \textbf{terminal payoff function} $F:Y \rightarrow \R$, a {\bf
running payoff function} $f:X \smallsetminus Y \rightarrow \R$,
and an \textbf{initial state} $x_0 \in X$.

Game play is as follows: a token is initially placed at position
$x_0$. At the $k\thh$ step of the game, a fair coin is tossed, and
the player who wins the toss may move the token to any $x_k$ for
which $(x_{k-1},x_k)$ is a directed edge in her transition graph.
The game ends the first time $x_k \in Y$, and player~\I's payoff
is $F(x_k) + \sum_{i=0}^{k-1} f(x_i)$.  Player~\I seeks to maximize
this payoff, and since the game is zero-sum, player~\II seeks
to minimize it.

We will use the term \textbf{tug-of-war (on the graph with edges
$E$)} to describe the game in which $E:=E_\I = E_\II$ (i.e., players
have identical move options) and $E$ is undirected (i.e., all moves
are reversible). Generally, our results pertain only to the
undirected setting. Occasionally, we will also mention some
counterexamples showing that the corresponding results do not hold
in the directed case.

In the most conventional version of tug-of-war on a graph, $Y$ is a
union of ``target sets'' $Y^\I$ and $Y^{\II}$, there is no running payoff ($f=0$) and $F$ is
identically $1$ on $Y^\I$, identically $0$ on $Y^{\II}$.  Players then try to ``tug'' the
game token to their respective targets (and away from their
opponent's targets) and the game ends when a target is reached.

A \textbf{strategy} for a player is a way of choosing the player's next
move as a function of all previously played moves and all previous coin tosses.
It is a map from the set of partially played games to moves (or in the case of a
\textbf{random strategy}, a probability distribution on moves).
Normally, one would
think of a good strategy as being Markovian, i.e., as a map from the
current state to the next move, but it is useful to allow more general
strategies that take into account the history.

Given two strategies $\mathcal S_{\I}, \mathcal S_{\II}$, let
$F_-(\mathcal S_{\I}, \mathcal S_{\II})$ and $F_+(\mathcal S_{\I},
\mathcal S_{\II})$ be the expected total payoff (including the
running payoffs received) at the termination of the game --- if the
game terminates with probability one and this expectation exists in
$[-\infty, \infty]$; otherwise, let $F_-(\mathcal S_{\I}, \mathcal
S_{\II})=-\infty$ and $F_+(\mathcal S_{\I}, \mathcal
S_{\II})=+\infty$.

The \textbf{value of the game for player~\I }is defined as
$\sup_{\mathcal S_\I} \inf_{\mathcal S_\II} F_-(\mathcal S_\I,
\mathcal S_\II)$. The \textbf{value for player~\II }is
$\inf_{\mathcal S_{\II}}\sup_{\mathcal S_{\I}} F_+(\mathcal
S_{\I},\mathcal S_{\II})$. We use the expressions $u_{\I}(x)$ and
$u_{\II}(x)$ to denote the values for players~\I and \II,
respectively, as a function of the starting state $x$ of the game.

Note that if  player~\I cannot force the game to end almost
surely, then $u_{\I}=-\infty$, and if  player~\II cannot force
the game to end almost surely, then $u_{\II}=\infty$. Clearly,
$u_{\I}(x) \leq u_{\II}(x)$. When $u_{\I}(x) = u_{\II}(x)$, we say
that the game has a \textbf{value}, given by
$u(x):=u_{\I}(x)=u_{\II}(x)$.

Our definition of value for player~\I penalizes player~\I severely
for not forcing the game to terminate with probability one, awarding
$-\infty$ in this case.

(As an alternative definition, one could define $F_-$, and hence
player~\I's value, by assigning payoffs to all of the
non-terminating sequences $x_0, x_1, x_2, \ldots$. If the payoff
function for the non-terminating games is a zero-sum
Borel-measurable function of the infinite sequence, then
 player~\I's value is equal to
 player \II's value in great generality \cite{MR1665779}; see also
\cite{MR2032421} for more on stochastic games.  The existence of a
value by our strong definition implies the existence and equality of
the values defined by these alternative definitions.)

Considering the two possibilities for the first coin toss yields
the following lemma, a variant of which appears in~\cite{MR1427981}.

\begin{lemma} \label{valuesareinfharmonic}
 The function $u=u_{\I}$ satisfies the equation
\begin{equation}\label{e.harm}
u(x) =  \frac12 \left( \sup_{y:(x,y)\in E_1} u(y) +
              \inf_{y:(x,y)\in E_2} u(y) \right) + f(x)
\end{equation}
for every non-terminal state $x \in X\setminus Y$ for which the
right-hand-side is well-defined, and $u_{\I}(x)=-\infty$ when the
right-hand-side is of the form $\frac12 (\infty +
(-\infty))+f(x)$. The analogous statement holds for $u_{\II}$,
except that $u_{\II}(x)=+\infty$ when the right-hand-side of
(\ref{e.harm}) is of the form $\frac12 (\infty + (-\infty))+f(x)$.
\end{lemma}

When $E=E_1=E_2$, the operator
$$\Delta_\infty u(x) :=
   \sup_{y:(x,y)\in E} u(y) +\inf_{y:(x,y)\in E} u(y) - 2u(x) $$
is called the \textbf{(discrete) infinity Laplacian}.  A
function $u$ is \textbf{infinity harmonic} if~\eqref{e.harm} holds and $f(x)=0$ at
all non-terminal $x \in X\setminus Y$. When $u$ is finite, this is
equivalent to $\Delta_{\infty}u=0$. However, it will be convenient
to adopt the convention that $u$ is infinity harmonic at $x$ if
$u(x)=+\infty$ (resp.\ $-\infty$) and the right hand side in~\eqref{e.harm} is also
$+\infty$ (resp.\ $-\infty$).  Similarly, it will be convenient to say ``$u$ is a solution to
$\Delta_\infty u = -2\,f$'' at $x$ if $u(x)=+\infty$ (resp.\ $-\infty$) and the right hand side in~\eqref{e.harm} is also
$+\infty$ (resp.\ $-\infty$).

In a tug-of-war game, it is natural to guess that the value
$u=u_{\I}=u_{\II}$ exists and is the \textit{unique\/} solution to
$$\Delta_\infty u(x) = -2f,$$ and also that (at least when $E$ is locally finite) player~\I's optimal
strategy will be to always move to the vertex that maximizes
$u(x)$, and player~\II's optimal strategy will be to always move
to the vertex that minimizes $u(x)$. This is easy to prove when
$E$ is undirected and finite and $f$ is everywhere positive or
everywhere negative. Subtleties arise in more general cases ($X$
infinite, $E$ directed, $F$ unbounded, $f$ having values of both
signs, etc.)

Our first theorem addresses the question of the existence of a value.

\begin{theorem} \label{towvalue}
A tug-of-war game with parameters $X, E, Y, F, f$ has a value whenever the following hold:
\begin{enumerate}
\item Either $f=0$ everywhere or $\inf f > 0$.
\item $\inf F >-\infty$.
\item $E$ is undirected.
\end{enumerate}
\end{theorem}
Counterexamples exist when any one of the three criteria is
removed. In Section~\ref{counterexamplesection} (a section devoted
to counterexamples) we give an example of a tug-of-war game
without a value, where  $E$ is undirected, $F=0$, and the running
payoff satisfies  $f > 0$ but $\inf f = 0$.

The case where $f=0$, $F$ is bounded, and $(X,E)$ is locally finite was proved
earlier in \cite{MR1685133}. That paper discusses an (essentially
non-random) game in which the two players bid for the
right to choose the next move.
That game, called the Richman game, has the same value as tug-of-war with $f=0$ where $F$ takes
the values $0$ and $1$.
Additionally, a simple and efficient algorithm for calculating the value
when $f=0$ and $(X,E)$ is finite is presented there.

\subsection{Tug-of-war on a metric space}

Consider now the special case that $(X,d)$ is a metric space,
$Y \subset X$, and Lipschitz functions $F:Y \rightarrow \R$ and
$f:X \smallsetminus Y \rightarrow \R$ are given. Let $E_\eps$ be
the edge-set in which $x \sim y$ if and only if $d(x,y) <\eps$ and
let $u^\eps$ be the value (if it exists) of the game played on $E_\eps$ with
terminal payoff $F$ and running payoff normalized to be $\eps^2 f$.

In other words, $u^\eps(x)$ is the value of the following
two-player zero-sum game, called \textbf{$\eps$-tug-of-war}: fix
$x_0=x\in X \smallsetminus Y$. At the $k\thh$ turn, the players
toss a coin and the winner chooses an $x_k$ with $d(x_k, x_{k-1})
<\eps$. The game ends when $x_k \in Y$, and player~\I's payoff is
$F(x_k) + \eps^2\sum_{i=0}^{k-1} f(x_i)$.

When the limit $u:=\lim_{\eps \rightarrow 0} u^\eps$ exists
pointwise, we call $u$ the \textbf{continuum value} (or just
``value'') of the quintuple $(X,d,Y,F,f)$.  We define the
\textbf{continuum value for player~\I }(or \II) analogously.

The reader may wonder why we have chosen not to put an edge in $E_\eps$
between $x$ and $y$ when $d(x,y)=\eps$ exactly. This
choice has some technical implications. Specifically,
we will compare the $\eps$-game with the $2\,\eps$-game.
If $x,z$ are such that $d(x,z)\le 2\,\eps$, then
in a length-space it does not follow that there
is a $y$ such that $d(x,y)\le\eps $ and $d(y,z)\le \eps$.
However, it does follow if you replace the weak inequalities
with strong inequalities throughout.

\medskip
We prove the following:

\begin{theorem} \label{continuumvalueexists}
  Suppose $X$ is a length space, $Y\subset X$ is non-empty, $F:Y
  \rightarrow \R$ is bounded below and has an
  extension to a uniformly continuous function on $X$, and either $f:X
  \smallsetminus Y \rightarrow \R$ satisfies $f=0$ or all three of the
  following hold: $\inf |f| > 0$, $f$ is uniformly continuous, and $X$
  has finite diameter.  Then the continuum value $u$ exists and is a uniformly
  continuous function extending $F$. Furthermore, $\|u-u^\eps\|_\infty\to 0$ as $\eps\searrow 0$.
  If $F$ is Lipschitz, then so is $u$.
  If $F$ and $f$ are Lipschitz, then $\|u-u^\eps\|_\infty = O(\eps)$.
\end{theorem}

The above condition that $F:Y\rightarrow\R$ extends to a uniformly continuous
function on $X$ is equivalent to having $F$ uniformly continuous on $Y$
and ``Lipschitz on large scales,'' as we prove in Lemma~\ref{extensionconditions} below.

We will see in Section~\ref{sub:pospay} that this fails in
general when $f > 0$ but $\inf f = 0$.
When $f$ assumes values
of both signs, it fails even when $X$ is a closed disk in
$\R^2$, $Y$ is its boundary and $F=0$.
In Section~\ref{sub:noval} we show by means of an example that in such circumstances
it may happen that $u_\I^\eps\ne u_\II^\eps$ and moreover,
$\liminf_{\eps\searrow 0} \|u_\I^\eps - u_\II^\eps\|_\infty>0$.

\subsection{Absolutely minimal Lipschitz extensions}
Given a metric space $(X,d)$, a subset  $Y \subset X$ and
a function $u:X \rightarrow \R$, we write $\Lip_Yu = \sup_{x,y \in Y}
|u(y)-u(x)|/d(x,y)$ and $\Lip\, u  = \Lip_Xu$.   Thus $u$ is
Lipschitz iff $\Lip\, u  < \infty$.  Given $F: Y \to \R$, we say
that $u:X\rightarrow \R$ is a \textbf{minimal extension of $F$} if
$\Lip_X u = \Lip_YF$ and $u(y) = F(y)$ for all $y \in Y$.

It is well known that for any metric space $X$, any Lipschitz $F$ on a
subset $Y$ of $X$ admits a minimal extension.  The largest and
smallest minimal extensions (introduced by McShane \cite{mcshane} and
Whitney \cite{MR1501735} in the 1930's) are respectively
$$\inf_{y \in Y} \left[F(y) + \Lip_Y F \,d(x,y)\right]
\ \ \ \ \ \text{and}\ \ \ \ \
\sup_{y \in Y} \left[F(y) - \Lip_Y F \,d(x,y)\right].$$

We say $u$ is an \textbf{absolutely minimal (AM) extension of $F$}
if $\Lip\, u <\infty$ and $\Lip_U u = \Lip_{\partial U} u$ for
every open set $U \subset X \smallsetminus Y$.
We say that $u$ is AM on $U$ if it is defined on
$\overline U$ and is an AM extension of its restriction to
$\partial U$.
AM extensions were
first introduced by Aronsson in 1967 \cite{MR0217665} and have
applications in engineering and image processing (see
\cite{MR2083637} for a recent survey).

We prove the following:

\begin{theorem} \label{uniqueAM}
  Let $X$ be a length space and let $F:Y\rightarrow \R$
be Lipschitz, where $\emptyset\ne Y\subset X$.
If $\inf F>-\infty$, then the continuum value
function $u$ described
  in Theorem~\ref{continuumvalueexists} (with $f=0$) is an AM
  extension of $F$.  If $F$ is also bounded, then $u$ is the
  \textit{unique\/} AM extension of $F$.
\end{theorem}

We present in the counterexample section, Section
\ref{counterexamplesection}, an example in which $F$ is Lipschitz,
non-negative, and unbounded, and although the continuum value is an AM extension, it is
not the only AM extension.

Prior to our work, the existence of AM extensions in the above
settings was known only for \textit{separable\/} length spaces
\cite{MR1884349} (see also \cite{MR1719573}).  The uniqueness in
Theorem~\ref{uniqueAM} was known only in the case that $X$ is the
closure of a bounded domain $U\subset \R^n$ and $Y=\partial U$.
(To deduce this case from Theorem~\ref{uniqueAM}, one needs to
replace $X$ by the smallest closed ball containing $\overline U$,
say.)
  Three uniqueness proofs in this setting have been
published, by Jensen \cite{MR1218686}, by Barles and Busca
\cite{MR1876420}, and by Aronsson, Crandall, and Juutinen
\cite{MR2083637}.  The third proof generalizes from the Euclidean norm
to uniformly convex norms.

Our proof applies to more general spaces because it invokes no outside
theorems from analysis (which assume existence of a local Euclidean
geometry, a measure, a notion of twice differentiability, etc.), and
relies only on the structure of $X$ as a length space.

As noted in \cite{MR2083637}, AM extensions do not generally
exist on metric spaces that are not length spaces. (For example, if $X$ is the
$L$-shaped region $\{0\} \times [0,1] \cup [0,1] \times \{0\}
\subset \R^2$, with the Euclidean metric, and
$Y = \{(0,1), (1,0) \}$, then no non-constant $F:Y\rightarrow
\R$ has an AM extension.
Indeed, suppose that $u:X\rightarrow\R$ is an AM extension of
$F:Y\rightarrow\R$. Let $a:=u(0,0)$, $b:=u(0,1)$ and
$c:= u(1,0)$. Then, considering $U=\{0\}\times (0,1)$,
it follows that $u(0,s)=a+s\,(b-a)$.
Likewise, $u(s,0)=a+s\,(c-a)$. Now taking
$U_\epsilon:=\{0\}\times[0,1)\cup [0,\epsilon)\times\{0\}$,
 we see that $\lim_{\eps\searrow 0} \Lip_{U_\eps}u=|b-a|$.
Hence $|c-a|\le |b-a|$. By symmetry, $|c-a|=|b-a|$.
Since $F$ is assumed to be nonconstant, $b\ne c$, and hence
$c-a=a-b$. Then
$\bigl|u(0,s)-u(s,0)\bigr|/\bigl(\sqrt 2 s\bigr) = \sqrt 2\,|b-a|$,
which contradicts $\lim_{\eps\searrow 0} \Lip_{U_\eps}u=|b-a|$.)

One property that makes length spaces special is that
the Lipschitz norm is determined locally. More precisely,
if $W\subset X$ is closed, then either $\Lip_W u=\Lip_{\partial W}u$
or for every $\delta>0$
$$
\sup\Bigl\{\frac{|u(x)-u(y)|}{d(x,y)}:
x,y\in W,\, 0<d(x,y)<\delta\Bigr\}= \Lip_Wu\,.
$$
  The definition of AM is inspired by the notion that if $u$ is the
``tautest possible'' Lipschitz extension of $F$, it should be
tautest possible on any open $V \subset X \smallsetminus Y$, given
the values of $u$ on $\partial V$ and ignoring the rest of the
metric space.  Without locality, the rest of the metric space
cannot be ignored (since long-distance effects may change the
global Lipschitz constant), and the definition of AM is less
natural. Another important property of length spaces is that the
graph distance metric on $E_\eps$ scaled by $\eps$ approximates the
original metric, namely, it is within $\eps$ of $d(\cdot, \cdot)$.

\subsection{Infinity Laplacian on $\R^n$} \label{inflapintro}

The continuum version of the \textbf{infinity Laplacian} is defined
for $C^2$ functions $u$ on domains $U \subset \R^n$ by
$$\Delta_\infty u = |\nabla u|^{-2} \sum_{i,j} u_{x_i} u_{x_ix_j} u_{x_j}.$$
This is the same as $\eta^T H \eta$, where $H$ is the Hessian of $u$ and
$\eta=\nabla u/|\nabla u|$.
Informally, $\Delta_\infty u$ is the second derivative of $u$ in the
direction of the gradient of $u$.  If $\nabla u(\vec x) = 0$ then, $\Delta_\infty u(\vec x)$ is undefined;
however, we adopt the convention that if the second derivative of $u(\vec x)$ happens to be
the same in \textit{every\/} direction (i.e, the matrix $\{u_{x_i x_j} \}$ is $\lambda$ times the
identity), then $\Delta_\infty u(\vec x) = \lambda$, which is the second derivative in any
direction.  (As mentioned above, some texts on infinity harmonic functions define $\Delta_\infty$ without the
normalizing factor $|\nabla u|^{-2}$.  When discussing viscosity solutions to
$\Delta_\infty u = 0$, the two definitions are equivalent.  The
fact that the normalized version is sometimes undefined when $\nabla u = 0$ turns out not
to matter because it \textit{is\/} always well-defined at $\vec x$ when $\phi$ is a {\bf
cone function}, i.e., it has the form $a|\vec x-\vec z|+b$ for $a, b \in \R$ and
$\vec z \in \R^n$ with $\vec z \neq \vec x$,
and viscosity solutions can be defined via comparison with cones, see Section~\ref{viscosityonlengthspace}.)
As in the discrete setting, $u$ is \textbf{infinity harmonic} if $\Delta_\infty u = 0$.

While discrete infinity-harmonic functions are a recent concept,
introduced in finite-difference schemes for approximating continuous
infinity harmonic functions \cite{oberman}, related notions
of value for stochastic games are of course much older.  The
continuous infinity Laplacian appeared first in the work of Aronsson 
\cite{MR0217665} and has been very thoroughly studied
\cite{MR2083637}.  Key motivations for studying this
operator are the following:
\begin{enumerate}
\item \textbf{AM extensions:}
Aronsson proved that $C^2$ extensions $u$ on domains $U\subset\R^n$
(of functions $F$ on $\partial U$) are infinity harmonic if and only if they are AM.
\item \textbf{$p$-harmonic functions:}
As noted by Aronsson~\cite{MR0217665}, the infinity Laplacian is the formal limit,
 as $p \rightarrow \infty$ of the (properly normalized)
$p$-Laplacians.  Recall that \textbf{$p$-harmonic} functions, i.e.\
minimizers $u$ of $\int |\nabla u(\vec x)|^p \,d\vec x$ subject to boundary
conditions, solve the \textbf{Euler-Lagrange equation}
$$\nabla\cdot(|\nabla u|^{p-2}\nabla u) = 0,$$ which can be rewritten
$$|\nabla u|^{p-2} \left(\Delta u + (p-2)\Delta_\infty u\right) =
0,$$ where $\Delta$ is the ordinary Laplacian.  Dividing by
$|\nabla u|^{p-2}$, we see that (at least when $|\nabla u| \neq 0$) $p$-harmonic functions satisfy
$\Delta_p u = 0$, where $\Delta_p := \Delta_\infty + (p-2)^{-1}
\Delta$; the second term vanishes in the large $p$ limit.
It is not too hard to see that as $p$ tends to infinity, the
Lipschitz norm of any limit of the $p$-harmonic functions extending $F$ will
be $\Lip_{\partial U}F$; so it is natural to guess (and was proved in
\cite{MR1155453}) that as $p$ tends to infinity the $p$-harmonic
extensions of $F$ converge to a limit that is both absolutely
minimal and a viscosity solution to $\Delta_\infty u= 0$.
\end{enumerate}

In the above setting, Aronsson also proved that there always exists an AM extension,
and that in the planar case $U\subset\R^2$, there exists at most
one $C^2$ infinity harmonic extension; however $C^2$
infinity harmonic extensions do not always exist \cite{MR0237962}.

To define the infinity Laplacian in the non-$C^2$ setting requires
us to consider weak solutions; the right notion here is that of
\textit{viscosity solution}, as introduced by Crandall and Lions (1983) \cite{MR690039}.
Start by observing that
if $u$ and $v$ are $C^2$ functions, $u(\vec x) = v(\vec x)$, and
$v \geq u$ in a neighborhood of $\vec x$, then $v-u$ has a
local minimum at $x$,
whence $\Delta_\infty v(\vec x) \geq
\Delta_\infty u(\vec x)$ (if both sides of this inequality
are defined).
 This comparison principle (which has analogs
for more general degenerate elliptic PDEs \cite{MR1462698}) suggests
that if $u$ is not $C^2$, in order to define $\Delta_\infty u(x)$ we want to
compare it to $C^2$ functions $\phi$ for which $\Delta_\infty \phi(x)$
is defined.
Let $S(x)$ be the set of real valued functions $\phi$
defined and $C^2$ in a neighborhood of $x$ for which $\Delta_\infty\phi(x)$
has been defined; that is, either $\nabla \phi(x)\ne 0$ or
$\nabla\phi(x)=0$ and the limit
$\Delta_\infty \phi(x):=
\lim_{x'\to x} 2\,\frac{\phi(x')-\phi(x)}{|x'-x|^2}$
exists.

\medskip

\noindent\textbf{Definition}.
Let $X$ be a domain in $\R^n$ and let $u: X \to \R$ be continuous.
Set
\begin{equation} \label{del+}
\Delta_{\infty}^{+}u(x) =\inf \{ \Delta_{\infty}\phi(x) \, :  \,
\phi\in S(x)\text{ and $x$ is a local minimum of }\phi-u\}\,.
\end{equation}
Thus $u$ satisfies $\Delta_{\infty}^{+}(u) \ge g$ in a domain $X$, iff
every $\phi \in C^2$ such that $\phi-u$ has a local minimum at some
$x\in X$ satisfies $\Delta^+_{\infty}\phi(x) \ge g(x)$.  In this case
$u$ is called a \textbf{viscosity subsolution} of
$\Delta_{\infty}(\cdot) =g$.
Note that if $\phi\in C^2$, then $\Delta_\infty^+\phi=\Delta_\infty\phi$
wherever $\nabla\phi\ne 0$.

Similarly, let
\begin{equation} \label{del-}
\Delta_{\infty}^{-}u(x) =\sup \{ \Delta_{\infty}\phi(x) \, :  \,
\phi\in S(x)\text{ and $x$ is a local maximum of }\phi-u\}\,.
\end{equation}
and call $u$  a  \textbf{viscosity supersolution} of $\Delta_{\infty}(\cdot) =g$
iff  $\Delta_{\infty}^{-}u \le g$ in $X$.

Finally, $u$ is a \textbf{viscosity solution}  of
$\Delta_{\infty}(\cdot) =g$ if $\Delta_{\infty}^{-}u \le g \le
\Delta_{\infty}^{+}u$ in $X$  (i.e., $u$ is both a supersolution and
a subsolution).

Here is a little caveat.
At present, we do not know how to show that $\Delta_\infty u=g$ in the viscosity
sense determines $g$. For example, if $u$ is Lipschitz, $g_1$ and $g_2$ are
continuous, and $\Delta_\infty u=g_j$ holds for $j=1,2$
(in the viscosity sense),  how does one prove that $g_1=g_2$?

\medskip

The following result of Jensen (alluded to above) is now well known
\cite{MR0217665, MR1218686, MR2083637}: If $X$ is a domain in $\R^n$
and $u:X \rightarrow \R$ continuous, then $\Lip_Uu = \Lip_{\partial U}u <
\infty$ for every bounded open set $U\subset\overline U \subset X$ (i.e., $u$ is AM)
if and only if $u$ is a viscosity solution to $\Delta_\infty u = 0$ in $X$.

\medskip

Let $A\subset Y\subset X$, where $A$ is closed, $Y\ne\emptyset$ and $X$ is a
length space. If $x\in X$, one can define
the \textbf{$\infty$-harmonic measure
of $A$ from $x$} as the infimum of $u(x)$ over all functions
$u:X\to [0,\infty)$ that are Lipschitz on $X$, AM in $X\setminus Y$ and satisfy
$u\ge 1$ on $A$.  This quantity will be denoted by
$\omega_\infty(A)=\omega^{(x,Y,X)}_\infty(A)$.
In Section~\ref{section:porous} we prove

\begin{theorem}\label{t.hmeasure}
Let $X$ be the unit ball in $\R^n$, $n>1$, let $Y=\partial X$, let $x$ be the center
of the ball,
and for each $\delta>0$ let $A_\delta\subset Y$ be a spherical cap of radius $\delta$
(of dimension $n-1$). Then
$$
c\,\delta^{1/3}\le \omega_\infty(A_\delta)\le C\,\delta^{1/3}\,,
$$
where $c,C>0$ are absolute constants (which do not depend on $n$).
\end{theorem}

Numerical calculations~\cite{oberman} had
suggested that in the setting of the theorem $\omega_\infty(A_\delta)$ tends to $0$ as $\delta\to 0$, but this was only
recently proved \cite{evans-yu}, and the proof did not yield any
quantitative information on the rate of decay.
In contrast to our other theorems in the paper, the proof of this theorem
does not use tug-of-war. The primary tool is comparison with a
specific AM function in $\R^2\setminus\{0\}$ with
decay $r^{-1/3}$ discovered by Aronsson~\cite{MR850366}.

\subsection{Quadratic comparison on length spaces} \label{viscosityonlengthspace}

To motivate the next definition, observe that for continuous functions
$u:\R \to \R$, the inequality $\Delta_\infty^+ u \ge 0$ reduces to convexity.
The definition of convexity requiring a function to lie below its chords
has an analog, \textit{comparison with cones}, which characterizes
infinity-harmonic functions.
Call the function $\phi(y)=b\,|y-z|+c$ a \textbf{cone} based at  $z \in \R^n$.
For an open $U \subset \R^n$, say that a continuous
$u: U \to \R$ satisfies \textbf{comparison with cones from above} on $U$ if for
every open
$W\subset\overline W  \subset  U$ for every $z\in \R^n\setminus W$, and
for every cone $\phi$ based at $z$
such that the inequality $u \le \phi$ holds  on $\partial W$,
the same inequality is valid throughout $W$.
\textbf{Comparison with cones from below} is defined similarly using the inequality $u \ge \phi$.

Jensen~\cite{MR1218686} proved that viscosity solutions to $\Delta_\infty u = 0$ for domains in $\R^n$
satisfy comparison with cones (from above and below), and Crandall, Evans, and
Gariepy~\cite{MR1861094} proved that a function on $\R^n$ is absolutely minimal in a bounded domain $U$
if and only if it satisfies comparison with cones in $U$.

Champion and De Pascale \cite{champdepascale} adapted this definition
to length spaces, where cones are replaced
by functions of the form $\phi(x)=b\,d(x,z)+c$ where $b>0$.
Their precise definition is as follows.
Let $U$ be an open subset of a length-space $X$ and let
$u: \overline U\to\R$ be continuous. Then $u$
is said to satisfy \textbf{comparison with distance functions from above} on $U$
if for every open $W\subset U$, for every $z\in X\setminus W$
for every $b\ge 0$ and for every $c\in \R$ if
$u(x)\le b\,d(x,z)+c$ holds on $\partial W$, then it also holds in $W$.
The function $u$ is said to satisfy \textbf{comparison with distance functions from below}
if $-u$ satisfies comparison with distance functions from above.
Finally, $u$ satisfies \textbf{comparison with distance functions} if it satisfies
comparison with distance functions from above and from below.

The following result from \cite{champdepascale} will be used in the proof of Theorem~\ref{uniqueAM}:

\begin{lemma}[\cite{champdepascale}] \label{AMiffCDF}
Let $U$ be an open subset of a length space.
  A continuous $u:\overline U\to \R$ satisfies comparison with distance functions in
  $U$ if and only if it is AM in $U$.
\end{lemma}

\bigskip
To study the inhomogenous equation $\Delta_\infty u=g$,
because $u$ will in general have non-zero second derivative (in its
gradient direction), it is natural to extend these definitions to
comparison with functions that have a quadratic term.

\medskip
\noindent\textbf{Definitions.} Let
 $Q(r)=a r^2+b r+c$ with $r,a,b,c\in\R$ and let $X$ be a length space.
\begin{itemize}
\item Let $z \in X$.
We call the function  $\phi(x) = Q(d(x,z))$ a \textbf{quadratic distance function}
(centered at $z$).
\item We say a quadratic distance function $\phi(x)=Q(d(x,z))$
is \textbf{\increasing} (in distance from $z$) on an open
set $V\subset X$ if either (1) $z\notin V$ and for every $x \in V$, we have $Q'(d(x,z))
> 0$, \textit{or\/} (2) $z\in V$ and $b=0$ and $a>0$. Similarly, we say
a quadratic distance function $\phi$ is \textbf{\decreasing} on $V$ if
$-\phi$ is \increasing on $V$.

\item
If $u:U\to\R$ is a continuous function defined on an
open set $U$ in a length space $X$,
we say that $u$
satisfies \textbf{$g$-quadratic comparison} on $U$
if the following two conditions hold:
\begin{enumerate}
\item \textsc{$g$-quadratic comparison from above}:
 For every open $V \subset\overline V\subset U$ and \increasing quadratic
distance function $\phi$ on $V$ with quadratic term $a \leq
\inf_{y \in V} \frac{g(y)}{2}$, the inequality $\phi \geq u$ on
$\partial V$ implies $\phi \geq u$ on $V$.
\item \textsc{$g$-quadratic comparison from below}:
 For every open $V \subset\overline V\subset U$ and \decreasing
quadratic distance function $\phi$ on $V$ with quadratic term $a
\geq \sup_{y \in V} \frac{g(y)}{2}$, the  inequality $\phi \leq u$
on $\partial V$ implies $\phi \leq u$ on $V$.
\end{enumerate}
\end{itemize}

The following theorem is proved in Section~\ref{visc-comp}.
\begin{theorem}  \label{thm-visc-comp}
  Let $u$ be a real-valued continuous function on a bounded domain $U$ in
  $\R^n$, and suppose that $g$ is a continuous function on $U$.  Then
  $u$ satisfies $g$-quadratic comparison on $U$ if and only if $u$ is a viscosity
  solution to $\Delta_\infty u = g$ in $U$.
\end{theorem}
This equivalence motivates the study of functions satisfying quadratic comparison.
Note that satisfying $\Delta_\infty u(x)=g(x)$ in the viscosity sense
depends only on the local behavior of $u$ near $x$.
We may use Theorem~\ref{thm-visc-comp} to extend the definition
of $\Delta_\infty$ to length spaces; saying that $\Delta_\infty u=g$
on an open subset $U$ of a length space if and only if every
$x\in U$ has a neighborhood $V\subset U$ on which
$u$ satisfies $g$-quadratic comparison.
We warn the reader, however, that even for length spaces $X$ contained within $\R$,
there can be solutions to $\Delta_\infty u=0$ that do not satisfy comparison with distance
functions, (or $0$-quadratic comparison, for that matter): for example, $X=U=(0,1)$ and $u(x)=x$.
The point here is that when we take $V\subset (0,1)$ and compare $u$ with some function
$\phi$ on $\partial V$, the ``appropriate'' notion of the boundary $\partial V$
is the boundary in $\R$, not in $X$.

The continuum value of the tug-of-war game sometimes gives a construction
of a function satisfying $g$-quadratic comparison. We prove:

\begin{theorem} \label{deltainftyuisf}
  Suppose $X$ is a length space, $Y\subset X$ is non-empty, $F:Y
  \rightarrow \R$ is uniformly continuous and bounded, and either $f:X
  \smallsetminus Y \rightarrow \R$ satisfies $f=0$ or all three of the
  following hold: $\inf |f| > 0$, $f$ is uniformly continuous, and $X$
  has finite diameter.
  Then the continuum value $u$ is the unique continuous function
  satisfying $(-2f)$-quadratic comparison on $X\smallsetminus Y$ and
  $u=F$ on $Y$.
  Moreover, if $\tilde u:X\to\R$ is continuous,
  satisfies $\tilde u\ge F$ on $Y$ and $(-2f)$-quadratic comparison from below on
  $X\setminus Y$, then $\tilde u\ge u$ throughout $X$.
\end{theorem}

Putting these last two theorems together, we obtain
\begin{cor}\label{c.linf}
  Suppose $U\subset\R^n$ is a bounded open set, $F:\partial U\to\R$ is
  uniformly continuous, and $f:U\to\R$ satisfies either $f=0$ or else
  $\inf f>0$ and $f$ is uniformly continuous.  Then there is a unique
  continuous function $u:\overline U\to\R$ that is a viscosity
  solution to $\Delta_\infty u = -2f$ on $U$ and satisfies $u=F$ on
  $\partial U$.  This unique solution is the continuum value
  of tug-of-war on $(\overline U,d,\partial U,F,f)$.
\end{cor}

It is easy to verify that $F$ indeed satisfies the assumptions in
Theorem~\ref{deltainftyuisf}.
In order to deduce the Corollary, we may take the length space $X$ as a
ball in $\R^n$ which contains $U$ and extend $F$ to $X\setminus U$,
say. Alternatively, we may consider $U$ with its intrinsic metric
and lift $F$ to the completion of $U$.

We present in Section~\ref{counterexamplesection}
an example showing that the corollary may fail if $f$ is permitted
to take values of both signs. Specifically, the example
describes two functions $u_1,u_2$ defined in the
closed unit disk in $\R^2$ and having boundary values
identically zero on the unit circle
such that with some Lipschitz function $g$ we have $\Delta_\infty u_j=g$
for $j=1,2$ (in the viscosity sense), while $u_1\ne u_2$.
\medskip

The plan of the paper is as follows.
Section~\ref{s.tug} discusses the discrete tug-of-war on graphs
and proves Theorem~\ref{towvalue},
and Section~\ref{mainproof} deals with tug-of-war on length spaces
and proves Theorems~\ref{continuumvalueexists}, \ref{uniqueAM} and~\ref{deltainftyuisf}.
Section~\ref{section:porous} is devoted to estimates of $\infty$-harmonic measure.
In Section~\ref{counterexamplesection}, we present a few counter-examples
showing that some of the assumptions in the theorems we prove are necessary.
Section~\ref{visc-comp} is devoted to the proof of Theorem~\ref{thm-visc-comp}.
Section~\ref{limitingtrajectory} presents some heuristic argument
describing what the limiting trajectories of some $\epsilon$-tug-of-war games on
domains in $\R^n$ may look like and states a question regarding the length of the game.
We conclude with additional open problems in Section~\ref{sec:open}.

\section{Discrete game value existence}\label{s.tug}

\subsection{Tug-of-war on graphs without running payoffs}

In this section, we will generally assume $E=E_{\I}=E_{\II}$ is undirected and
connected and $f=0$.

\medskip

Though we will not use this fact, it is interesting to
point out that
in the case where $E$ is finite, there is a simple
algorithm from~\cite{MR1685133} which calculates the value $u$ of the game
and proceeds as follows.
Assuming the value $u(v)$ is already calculated at some
set $V'\supset Y$ of vertices, find a path $v_0,v_1,\dots,v_k$, $k>1$,
with interior vertices $v_1,\dots,v_{k-1}\in X\setminus V'$ and endpoints $v_0,v_k\in V'$
which maximizes $\bigl(u(v_k)-u(v_0)\bigr)/k$
and set $u(v_i)=u(v_0)+i\,\bigl(u(v_k)-u(v_0)\bigr)/k$ for $i=1,2,\dots,k-1$.
Repeat as long as $V'\ne X$.

\medskip

Recall that $u_\I$
is the value function for player~\I.

\begin{lemma} \label{smallestinfharmonic}
Suppose that $\inf_{Y} F
 >-\infty$ and $f=0$.
 Then $u_\I$ is the smallest $\infty$-harmonic function bounded
 below on $X$ that extends $F$. More generally, if $v$ is an
 $\infty$-harmonic function which is bounded from below on $X$
 and $v \ge F$ on $Y$, then $v \ge u_\I$ on $X$.
\end{lemma}
Similarly, if $F$ is bounded above on $Y$, then $u_\II$ is the
largest $\infty$-harmonic function bounded
 above on $X$ that extends $F$.
\begin{proof}
  Player~\I could always try to move closer to some specific point $y\in Y$.
Since in almost
every infinite sequence of fair coin tosses there will be a
time when the number of tails exceeds the number of heads by $d(x_0,y)$,
this ensures that the game terminates a.s.,
and we have $u_\I \geq \inf_Y F > -\infty$.
Suppose that $v\ge F$ on $Y$ and $v$ is $\infty$-harmonic on
$X$. Given $\delta>0$, consider an arbitrary strategy for player~\I
and let \II play a strategy that at step $k$ (if \II wins
the coin toss), selects a state where $v(\cdot)$ is within $\delta
2^{-k}$ of its infimum among the states \II could move to. We
will show that the expected payoff for player~\I is at most
$v(x_0)+\delta$.    We may assume the
game terminates a.s.\ at a time $\tau<\infty$. Let $\{x_j\}_{j \ge
0}$ denote the random sequence of states encountered in the game.
 Since $v$ is $\infty$-harmonic, the sequence $M_k=v(x_{k\wedge
\tau}) + \delta 2^{-k}$ is a supermartingale. Optional sampling
and Fatou's lemma imply that $v(x_0)+\delta =M_0 \ge \E[M_\tau] \ge
\E[F(x_\tau)]$. Thus $u_\I\le v$.  By Lemma~\ref{valuesareinfharmonic}, this
completes the proof.
\end{proof}

Now, we prove the first part of Theorem~\ref{towvalue}: When the graph
$(X,E)$ is locally finite and $Y$ is finite, this was proven in
\cite[Thm.~17]{MR1685133}.

\begin{theorem} \label{uniqueinfharmongraph} Suppose that $(X,E)$ is
connected, and $Y\ne\emptyset$.  If $F$ is bounded
  below (or above) on $Y$ and $f=0$, then $u_\I=u_\II$, so the game has a value.
\end{theorem}

Before we prove this, we discuss several counterexamples that
occur when the conditions of the theorem are not met.  First,
this theorem can fail if $E$ is directed.  A trivial
counterexample is when $X$ is finite and there is no directed path from
the initial state to $Y$.

If $X$ is infinite and $E$ is directed, then there are
counterexamples to Theorem~\ref{uniqueinfharmongraph} even when
every vertex lies on a directed path towards a terminal state. For
example, suppose $X=\N$ and $Y=\{0\}$, with $F(0)=0$.  If $E$
consists of directed edges of the form $(n,n-1)$ and $(n,n+2)$,
then \II may play so that with positive probability the game
never terminates, and hence the value for
 player~\I is by definition $-\infty$.

\label{s.z2} Even in the undirected case, a game may not have a
value if $F$ is not bounded either from above or below.  The
reader may check that if $X$ is the integer lattice $\Z^2$ and the
terminal states are the $x$-axis with $F((x,0)) = x$, then the
players' value functions are given by $u_{\I}((x,y)) = x-|y|$ and
$u_{\II}((x,y)) = x +|y|$.  (Roughly speaking, this is because, in
order to force the game to end in a finite amount of time, player \I
has to ``give up'' $|y|$ opportunities to move to the right.)
Observe also that in this case any linear function which agrees
with $F$ on the $x$-axis is $\infty$-harmonic.

\begin{figure}
\centerline{\includegraphics[height=1.5in]{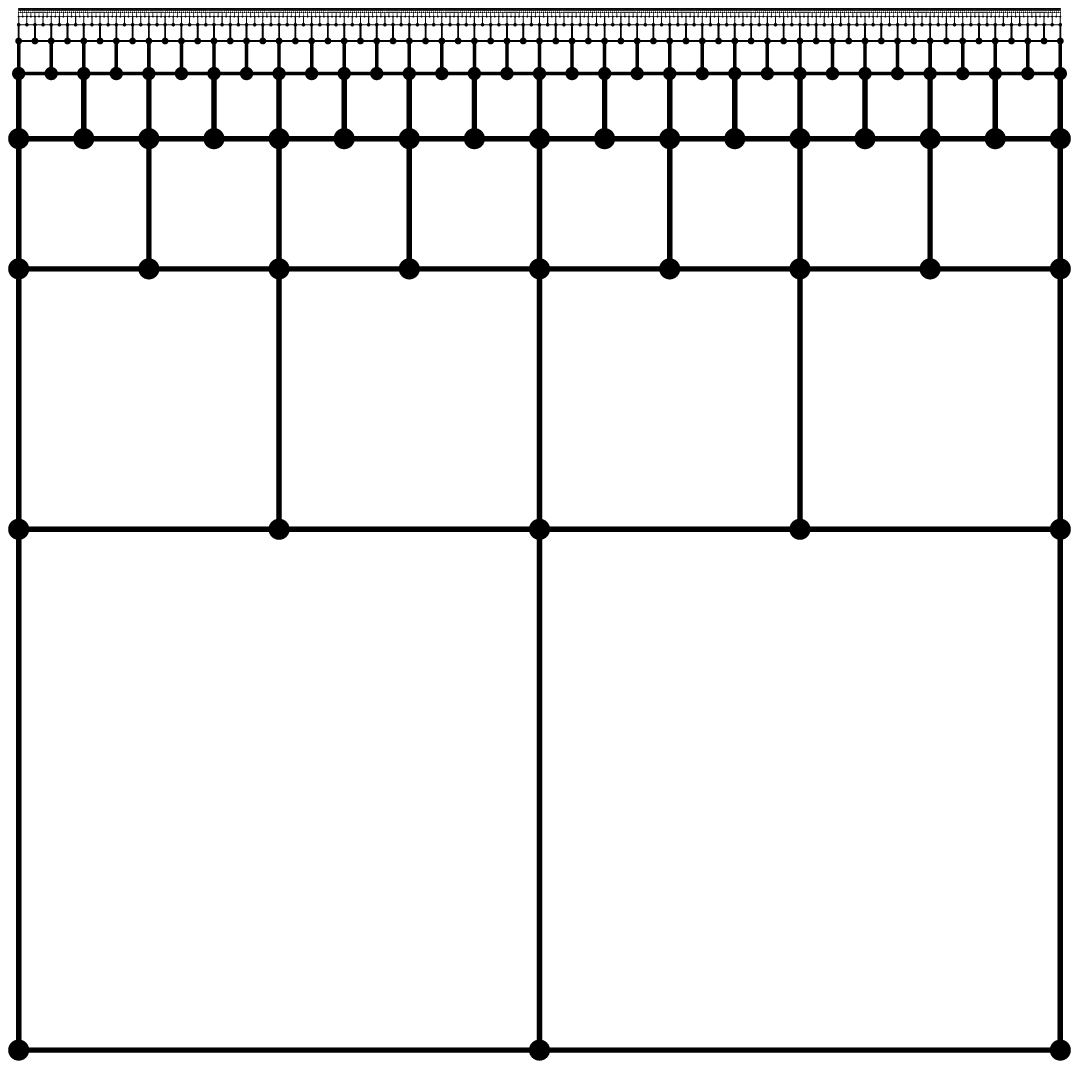}}
\caption
{A graph for which the obvious tug-of-war strategy of
maximizing/minimizing $u$ does poorly.\label{fig:pullup}}
\end{figure}

As a
final remark before proving this theorem, let us consider the
obvious strategy for the two players, namely for player~\I to always
maximize $u$ on her turn and for \II to always minimize $u$ on
her turn.  Even when $E$ is (undirected and) locally finite and the
payoff function satisfies $0\leq F\leq 1$, these obvious
strategies need not be optimal.  Consider, e.g., the game shown in
Figure~\ref{fig:pullup}, where $X\subset\R^2$ is given by
$X = \bigl\{v(k,j) :j=1,2,\dots; k
=0,1,\dots,2^j\bigr\}$, $v(k,j) =(2^{-j}k,1-2^{-j+1})$,
and $E$ consists of edges of the form $\{v(k,j),v(k+1,j)\}$ and
$\{v(k,j),v(2k,j+1)\}$.  The terminal states are on the left and
right edges of the square, and the payoff is $1$ on the left and
$0$ on the right.  Clearly the function  $u(v(k,j))=1-k/2^j$ (the
Euclidean distance from the right edge of the square) is
infinity-harmonic, and by Corollary~\ref{cor:uniqueinfharmongraph} below,
we have $u_\I=u=u_\II$.  The obvious strategy for player~\I is to
always move left, and for \II it is to always move right.
Suppose however that  player~\I
always pulls left and player~\II always pulls up.
It is easy to check that the probability that the game
ever terminates when starting from $v(k,j)$ is at most
$2/(k+2)$ (this function is a supermartingale under the
corresponding Markov chain).
Therefore the game continues forever with positive probability,
resulting in a payoff of $-\infty$ to player~\I.
Thus, a near-optimal strategy for player~\I must be able
to force the game to end, and must be prepared to make moves which
do not maximize $u$.
(This is a well known phenomenon, not particular to tug-of-war.)

\begin{proof}[Proof of Theorem~\ref{uniqueinfharmongraph}]
  If $F$ is bounded above but not below, then we may exchange the
  roles of players~$\I$ and $\II$ and negate $F$ to reduce to the
  case where $F$  is bounded below.
  Since $u_{\I} \leq u_{\II}$  always holds,
  we just need to show that $u_{\II}\le u_\I$.
  Since player~\I could always pull towards a point in $Y$ and thereby
  ensure that the game terminates, we have $u_\I \geq \inf_Y F > -\infty$.
 Let $u=u_\I$ and
write $\delta(x) = \sup_{y: y \sim x} |u(y)-u(x)|$. Let $x_0,x_1,\dots$ be the
sequence of positions of the game.
  For ease of exposition, we
  begin by assuming that $E$ is locally finite (so that the suprema
  and infima in the definition of the $\infty$-Laplacian definition
  are achieved) and that $\delta(x_0)>0$; later we will remove these assumptions.

To motivate the following argument, we make a few observations.
In order to prove that $u_\II\le u_\I$, we need to show that
player~\II can guarantee that the game terminates while also making sure
that the expected payoff is not much larger than $u(x_0)$.
These are two different goals, and it is not \textit{a priori\/} clear how to combine them.
To resolve this difficulty, observe that
$\delta(x_j)$ is non-decreasing in $j$ if players~\I and \II adopt the
strategies of maximizing (respectively, minimizing) $u$ at every
step.
As we will later see, this implies that the game terminates a.s.\ 
under these strategies. On the other hand, if
player~\I deviates from this strategy and thereby reduces $\delta(x_j)$,
then perhaps player~\II can spend some turns playing suboptimally
with respect to $u$ in order to increase $\delta$.
Let $X_0:=\{x\in X:\delta(x)\ge\delta(x_0)\}\cup Y$. For
$n=0,1,2,\dots$ let $j_n=\max\{j\le n:x_j\in X_0\}$ and
$v_n=x_{j_n}$, which is the last position in $X_0$ up to time $n$.
We will shortly describe a strategy for \II based on the idea
of backtracking to $X_0$ when not in $X_0$. If $v_n\ne x_n$, we
may define a backtracking move from $x_n$ as any move to a
neighbor $y_n$ of $x_n$ that is closer to $v_n$ than $x_n$ in the
subgraph $G_n\subset (X,E)$ spanned by the vertices
$x_{j_n},x_{j_n+1},\dots,x_n$.
Here, ``closer'' refers to the graph metric of $G_n$. When \II plays the backtracking
strategy, she backtracks whenever not in $X_0$ and plays to a
neighbor minimizing $u_\I$ when in $X_0$.

Now consider the game evolving under any strategy for player~\I and
the backtracking strategy for \II. Let $d_n$ be the distance
from $x_n$ to $v_n$ in the subgraph $G_n$. Set
$$
m_n:=u(v_n)+\delta(x_0)\,d_n\,.
$$
It is clear that $u(x_n)\le m_n$, because there is a path of
length $d_n$ from $x_n$ to $v_n$ in $G_n$ and the change in $u$ across
any edge in this path is less than $\delta(x_0)$, by the
definition of $X_0$. It is easy to verify that $m_n$ is a
supermartingale, as follows. If $x_n\in X_0$, and player~\I plays,
then $m_{n+1}\le u(x_n)+\delta(x_n)= m_n+\delta(x_n)$, while if
 \II gets the turn then $m_{n+1}=u(x_n)-\delta(x_n)$. If
$x_n\notin X_0$ and \II plays, then $m_{n+1}=m_n-\delta(x_0)$.
If $x_n\notin X_0$ and player~\I plays not into $X_0$, then
$m_{n+1}\le m_n+\delta(x_0)$. The last case to consider is that
$x_n\notin X_0$ and player~\I plays into a vertex in $X_0$. In such a
situation,
$$
m_{n+1}=u(x_{n+1})\le u(x_n)+\delta(x_n)\le m_n+\delta(x_0)\,.
$$
Thus, indeed, $m_n$ is a supermartingale (bounded below).  Let $\tau$
denote the first time a terminal state is reached (so
$\tau=\infty$ if the game does not terminate). By the martingale
convergence theorem, the limit $\lim_{n\to\infty}m_{n \wedge
\tau}$ exists. But when player \II plays we have $m_{n+1}\le
m_n-\delta(x_0)$. Therefore,  the game must terminate with
probability $1$. The expected outcome of the game thus played is
at most $m_0=u(x_0)$. Consequently, $u_{\II}\le u$, which
completes the proof in the case where $E$ is locally finite and
$\delta(x_0)>0$.

Next, what if $E$ is not locally finite, so that suprema and
infima might not be achieved?  In this case, we fix a small
$\eta > 0$, and use the same strategy as above except that if
$x_n\in X_0$ and \II gets the turn, she moves to a neighbor at
which $u$ is at most $\eta 2^{-n-1}$ larger than its infimum value
among neighbors of $x_n$.  In this case, $m_n + \eta\,2^{-n}$ is
a supermartingale, and hence the expected payoff is at most
$u(x_0) + \eta$.  Since this can be done for any $\eta>0$, we
again have that $u_{\II} \le u $.

Finally, suppose that $\delta(x_0)=0$.
Let $y\in Y$, and
let player \II pull toward $y$ until the first time a vertex $x_0^*$ with $\delta(x_0^*)>0$
or $x_0^*\in Y$ is reached.
After that, \II continues as above.  Since $u(x_0)=u(x_0^*)$, this completes the proof.
\end{proof}

\begin{cor} \label{cor:uniqueinfharmongraph}
  If $E$ is connected, $Y\ne\emptyset$ and $\sup |F|<\infty$,
  then $u=u_\I=u_\II$ is the unique bounded $\infty$-harmonic function
  agreeing with $F$ on $Y$.
\end{cor}
\begin{proof} This is an immediate consequence of Lemma~\ref{smallestinfharmonic},
the remark that follows it, and
Theorem~\ref{uniqueinfharmongraph}.
\end{proof}

If $E=\{(n,n+1):n=0,1,2,\dots\}$, $Y=\{0\}$ and $F(0)=0$,
then $\tilde u(n)=n$ is an example of an (unbounded) $\infty$-harmonic
function that is different from $u$.

\subsection{Tug-of-war on graphs with running payoffs} \label{s.run} Suppose now
that $f \neq 0$.  Then the analog of Theorem
\ref{uniqueinfharmongraph} does not hold without additional
assumptions.  For a simple counterexample, suppose that $E$ is a
triangle with self loops at its vertices (i.e., a player may opt
to remain in the same position), that the vertex $v_0$ is a
terminal vertex with final payoff $F(v_0)=0$, and the running
payoff is given by $f(v_1)=-1$ and $f(v_2)=1$. Then the function
given by $u(v_0)=0$, $u(v_1)=a-1$, $u(v_2) = a+1$ is a solution to
$\Delta_{\infty} u = -2f$ provided $-1 \leq a \leq 1$. The reader
may check that $u_{\I}$ is the smallest of these functions and
$u_{\II}$ is the largest.  The gap of $2$ between $u_{\I}$ and
$u_{\II}$ appears because a player would have to give up a move
(sacrificing one) in order to force the game to end. This is
analogous to the $\Z^2$ example given in \S\ref{s.z2}.
 Both players are earning payoffs in the
interior of the game, and moving to a terminal vertex costs a player
a turn.

One way around this is to assume that $f$ is either uniformly positive
or uniformly negative, as in the following analog of Theorem
\ref{uniqueinfharmongraph}.  We now prove the second half of
Theorem~\ref{towvalue}:

\begin{theorem} \label{uniquesolutionongraph} Suppose that $E$ is
connected and $Y\ne\emptyset$.
 Assume that $F$ is bounded below and $\inf f >0$.
 Then $u_\I = u_\II$.
 If, additionally, $f$ and $F$ are bounded above,
 then any bounded solution $\tilde u$  to $\Delta_\infty \tilde u = -2f$
with the given
 boundary conditions is equal to $u$.
\end{theorem}

\begin{proof}
  By considering a strategy for player~\I that always pulls toward a specific
  terminal state $y$, we see that $\inf_X u_\I \ge \inf_{Y} F$.
By Lemma~\ref{valuesareinfharmonic},
$\Delta_\infty u_\I=-2f$ on $X\smallsetminus Y$.
 Let $\tilde u$ be any solution to $\Delta_\infty \tilde u=-2f$ on
$X\smallsetminus Y$
  that is bounded below on $X$ and has the given boundary values on $Y$.
\smallskip

\noindent\textbf{Claim:} $u_\II  \le \tilde u$.
In proving this, we may assume without loss of generality that
$G\setminus Y$ is connected, where $G$ is the graph $(X,E)$.
Then if $\tilde u=\infty$ at some vertex
in $X\setminus Y$, we also have $\tilde u=\infty$ throughout $X\setminus Y$,
in which case the claim is obvious.
Thus, assume that $\tilde u$ is finite on $X\setminus Y$.
Fix $\delta\in(0,\inf f)$ and let \II use the strategy that at
  step $k$, if the current state is $x_{k-1}$ and \II wins the coin
  toss, selects a state $x_k$ with $\tilde u(x_k)<\inf_{z:z \sim x_{k-1}}
  \tilde u(z)+2^{-k}\delta.$ Then for any strategy chosen by player~\I, the
  sequence $M_k=\tilde u(x_k)+2^{-k}\delta+\sum_{j=0}^{k-1} f(x_j)$ is a
  supermartingale bounded from below, which must converge a.s.\ to a finite
  limit.  Since $\inf f>0$, this also forces the game to terminate
  a.s.  Let $\tau$ denote the termination time.  Then $$
  \tilde u(x_0)+\delta=M_0 \ge \E(M_\tau) \ge
  \E\Bigl(\tilde u(x_\tau)+\sum_{j=0}^{\tau-1} f(x_j) \Bigr) \,.  $$
  Thus
  this strategy for \II shows that $u_\II(x_0) \le \tilde u(x_0)+\delta$.
  Since $\delta>0$ is arbitrary, this verifies the claim.
In particular, $u_\II \le u_I$ in $X$, so $u_\II=u_\I$.

  Now suppose that $\sup F<\infty$ and $\sup f < \infty$, and
  $\tilde u$ is a \textit{bounded\/} solution to $\Delta_\infty \tilde u=-2f$
  with the given boundary values. By the claim above,
  $\tilde u\ge u_{\II}=u_\I$. On the other
  hand, player~\I can play to maximize or nearly maximize $\tilde u$ in every
  move. Under such a strategy, he guarantees that by turn $k$ the
  expected payoff is at least $\tilde u(x_0)-\E[Q(k)]-\eps$, where
  $Q(k)$ is $0$ if the game has terminated by time $k$ and $\tilde
  u(x_k)$ otherwise. If the expected number of moves played is
  infinite, the expected payoff is infinite.  Otherwise,
  $\lim_k\E[Q(k)]=0$, since $\tilde u$ is bounded. Thus, $u_\II=u_{\I}\ge
  \tilde u$ in any case.
\end{proof}

\section{Continuum value of tug-of-war on a length space} \label{mainproof}

\subsection{Preliminaries and outline}

In this section, we will prove Theorem~\ref{continuumvalueexists},
Theorem~\ref{deltainftyuisf} and Theorem~\ref{uniqueAM}.  Throughout
this section, we assume $(X,d,Y,F,f)$ denotes a tug-of-war game, i.e.,
$X$ is a length space with distance function $d$, $Y\subset X$ is a
nonempty set of terminal states, $F:Y\to\R$ is the final payoff
function, and $f:X\smallsetminus Y\to\R$ is the running payoff
function.  We let $x_k$ denote the game state at time $k$ in
$\eps$-tug-of-war.

It is natural to ask for a continuous-time version of tug-of-war on a
length space.  Precisely and rigorously defining such a game (which
would presumably involve replacing coin tosses with white noise,
making sense of what a continuum no-look-ahead strategy means, etc.)\
is a technical challenge we will not undertake in this paper (though
we include some discussion of the small-$\eps$ limiting trajectory of
$\eps$ tug-of-war in the finite-dimensional Euclidean case in
Section~\ref{limitingtrajectory}).  But we can make sense of the
continuum game's \textit{value function\/} $u^0$ by showing that the
value function $u^\eps$ for the $\eps$-step tug-of-war game converges
as $\eps\to 0$.

The value functions $u^\eps$ do not satisfy any nice monotonicity
properties as $\eps\to 0$.  In the next subsection we define two
modified versions of tug-of-war whose values closely approximate
$u^\eps$, and which do satisfy a monotonicity property along sequences
of the form $\eps 2^{-n}$, allowing us to conclude that $\lim_n u^{\eps
  2^{-n}}$ exists.  Then we show that any such limit is a bounded below viscosity
solution to $\Delta_\infty u=-2f$, and that any
viscosity solution bounded below is an upper bound on such a limit, so
that any two such limits must be equal, which will allow us to prove
that the continuum limit $u^0=\lim_\eps u^\eps$ exists.

Because the players can move the game state almost as far as $\eps$,
either player can ensure that $d(x_k,y)$ is ``almost a
supermartingale'' up until the time that $x_k=y$.  When doing
calculations it is more convenient to instead work with a related
metric $d^\eps$ defined by
\begin{align*}
d^\eps(x,y) : = {}&\eps\times(\text{min \# steps from $x$ to $y$ using steps of length $<\eps$})\\
 ={}&\begin{cases}0&x=y\\ \eps+\eps \lfloor d(x,y)/\eps \rfloor & x\neq y\end{cases}
\end{align*}
Since $d^\eps$ is the graph distance scaled by $\eps$, it is in fact a
metric, and either player may choose to make $d^\eps(x_k,y)$ a
supermartingale up until the time $x_k=y$.

\subsection{\II-favored tug-of-war and dyadic limits}

We define a game called \textbf{II-favored $\eps$-tug-of-war}
that is designed to give a lower bound on player~\I's expected payoff.
It is related to ordinary $\eps$-step tug-of-war, but \II is given
additional options, and player~\I's running payoffs are slightly smaller.  At
the $(i+1)$\st step, player~\I chooses a point $z$ in $B_{\eps} (x_i)$ and a
coin is tossed. If player~\I wins the coin toss, the game position moves to a
point, of player \II's choice, in $(B_{2\eps}(z) \cap Y) \cup \{z\}$. (If
$d(z,Y) \ge 2\eps$, this means simply moving to $z$.)  If \II wins,
then the game position moves to a point in $B_{2\eps}(z)$
of \II's choice.  The game ends at the first time $\tau$ for which
$x_\tau \in Y$. Player~\I's payoff is then $-\infty$ if the game never
terminates, and otherwise it is
\begin{equation} \label{epsilonpayoffequation}
\eps^2\sum_{i=1}^{\tau} \inf_{y \in B_{2\eps}(z_i)} f(y) + F(x_\tau)
\end{equation}
where $z_i$ is the point that player~\I targets on the $i$\thh turn,
and $f(y)$ is defined to be zero if $y \in Y$.  We let $v^\eps$ be
the value for player~\I for this game.  Given a strategy for
player~\II in the ordinary $\eps$-game, player~\II can easily mimic
this strategy in the \II-favored $\eps$-game and do at least as
well, so $v^\eps\leq u_\I^\eps$.

Let $w^\eps$ be the value for player~\II of \I-favored
$\eps$-tug-of-war, defined analogously but with the roles
of player~\I and player~\II reversed (i.e., at each move, \II selects
the target less than $\eps$ units away, instead of player~\I, etc., and the
$\inf$ in the running payoff term in
Equation~\eqref{epsilonpayoffequation} is replaced with a $\sup$, and
games that never terminate have payoff $+\infty$).
For any $\eps>0$ we have $$v^\eps \leq u_\I^\eps \leq u_\II^\eps \leq w^\eps.$$

\begin{lemma} \label{epsilongameinequalities}
  For any $\eps > 0$,
  $$v^{2\eps} \leq v^\eps \leq u_\I^\eps \leq u_\II^\eps \leq w^\eps \leq w^{2 \eps}.$$
\end{lemma}

\begin{proof}
  We have already noted that $v^\eps \leq u_\I^\eps \leq u_\II^\eps \leq w^\eps$.  We
  will prove that $v^{2\eps}\leq v^\eps$; the inequality $w^\eps \leq
  w^{2 \eps}$ follows by symmetry.  Consider a strategy $S_\I^{2\eps}$
  for \II-favored $2\eps$-tug-of-war.  We define a strategy $S_\I^\eps$
  for player~\I for the \II-favored $\eps$ game that mimics
  $S_\I^{2\eps}$ as follows.  Whenever strategy $S_\I^{2\eps}$ would
  choose a target point $z$, player~\I ``aims'' for $z$ for one
  ``round,'' which we define to be the time until one of the players
  has won the coin toss two more times than the other player.  By
  ``aiming for $z$'' we mean that player~\I picks a target point that,
  in the metric $d^\eps$, is $\eps$ units closer to $z$ than the
  current point.  With probability $1/2$, player~\I gets two surplus
  moves before \II, and then the game position reaches $z$ (or a point
  in $Y \cap B_{4\eps}(z)$) before the game position exits
  $B_{4\eps}(z)$.  If player~\II gets two surplus moves before player~\I,
  then the game position will be in $B_{4\eps}(z)$.
  (See Figure~\ref{fig:tw}.)
\begin{figure}[h]
\centerline{\includegraphics{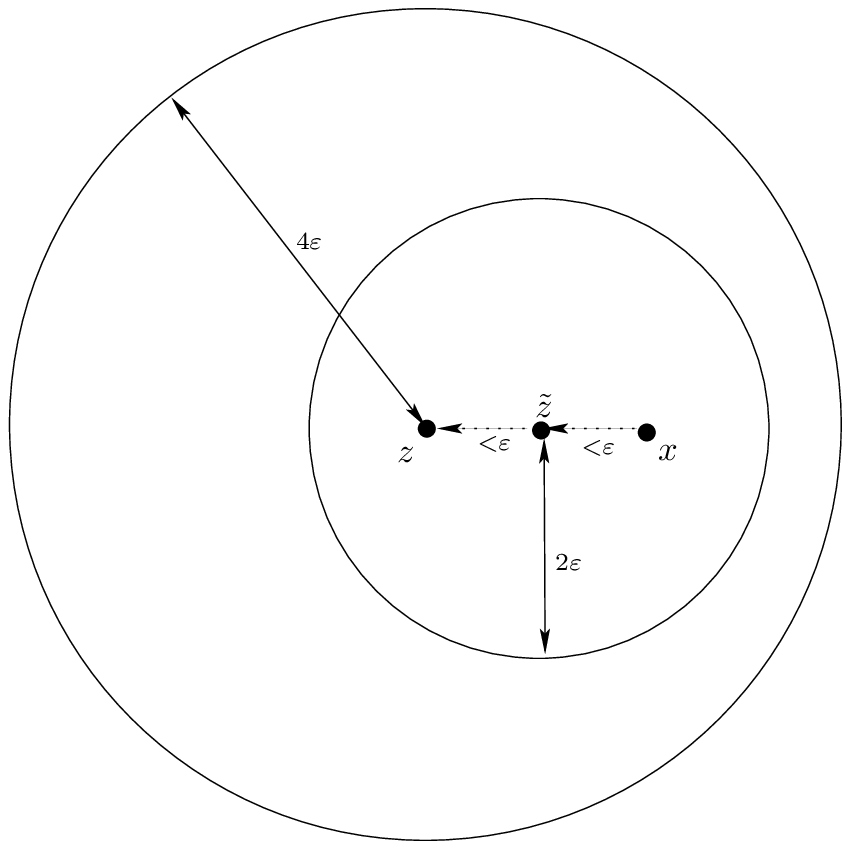}}
\caption{At the end of the round, player~\I reaches the target $z$ (or a point in $Y \cap B_{4\eps}(z)$) with probability at least $1/2$,
 and otherwise the state still remains within $B_{4\eps}(z)$.  During the first step of the round, when player~\I has target $\tilde z$, the running payoff is the infimum of $f$ over $B_{2\eps}(\tilde z) \subset B_{4\eps}(z)$, and during any step of the round the infimum is over a subset of $B_{4\eps}(z)$.  \label{fig:tw}}
\end{figure}
  The expected
  number of moves in this round of the \II-favored $\eps$ game is $4$;
  and the running payoff at each move is an $\eps^2$ times the infimum
  over a ball of radius $2\eps$ that is a subset of $B_{4\eps}(z)$ (as
  opposed to $(2\eps)^2$ times the infimum over the whole ball).
  Hence strategy $S_\I^\eps$ guarantees for the \II-favored $\eps$ game
  an expected total payoff that is at least as large as what
  $S_\I^{2\eps}$ guarantees for the \II-favored $2\eps$ game.
\end{proof}

Thus $v^\eps$ converges along dyadic sequences:
$v^{\eps/2^\infty}:=\lim_{n\to\infty} v^{\eps 2^{-n}}$ exists.  \textit{A
  priori\/} the subsequential limit could depend upon the choice of
the dyadic sequence, i.e., the initial $\eps$.

The same argument can be used to show that $v^{k\eps}\le v^\eps$ for positive integers $k$.

\subsection{Comparing favored and ordinary tug-of-war}

We continue with a preliminary bound on how far apart $v^\eps$ and
$w^\eps$ can be.
Let $\Lip^\eps_YF$ denote the Lipschitz constant of $F$ with respect to the
restriction of the metric $d^\eps$ to $Y$.
Since $d\le d^\eps$, the Lipschitz constant $\Lip_Y F$ of $F$
with respect to $d$ upper bounds $\Lip^\eps_YF$.

\begin{lemma}
  \label{vwgap}
  Let $\eps>0$.
  Suppose that $\Lip^\eps_YF<\infty$ and either
      \begin{enumerate} \item $f=0$ everywhere,\, or  \item $|f|$ is bounded above and $X$ has finite diameter.
      \end{enumerate}
  Then for each $x\in X$ and $y\in Y$,
  $$ v^\eps(x) \geq F(y) - 2\,\eps\, \Lip_Y^\eps F -
\bigl(\Lip_Y^\eps F + 2\, (\eps +\diam X)\sup|f|\bigr)\, d^\eps(x,y).$$
Such an expected payoff is guaranteed for player~\I if he adopts a
pull towards $y$ strategy, which at each move attempts to reduce $d^\eps(x_t,y)$.
Similarly a pull towards $y$ strategy for~\II gives,
  $$ w^\eps(x) \leq F(y) + 2\,\eps\, \Lip_Y^\eps F + \bigl(\Lip_Y^\eps F  +
2\, (\eps +\diam X)\sup|f|\bigr)\, d^\eps(x,y).$$
\end{lemma}
\begin{proof}
Let player~\I use a pull towards $y$ strategy.
  Let $\tau$ be the time at which $Y$ is reached, which will be finite
  a.s.  The distance $d^\eps(x_k,y)$ is a supermartingale, except
  possibly at the last step, where player~\II may have moved the game
  state to a terminal point up to a distance of $2\eps$ from the target, even
  if player~\I wins the coin toss.  Thus $\E[d^\eps(x_\tau,y)]< d^\eps(x,y) + 2\eps$,
whence $\E[F(x_\tau)] \geq F(y) - \bigl(d^\eps(x,y) + 2\eps\bigr)\,\Lip_Y^\eps F $.
If $f=0$, then this implies $v^\eps(x) \geq F(y) - (d^\eps(x,y) + 2\eps)\Lip_Y^\eps F $.
If $f \neq 0$ and $X$ has finite
  diameter, then the expected number of steps before the game
  terminates is at most the expected time that a simple random walk on
  the interval of integers $[0, 1+\lfloor \diam(X)/\eps \rfloor]$
 (with a self-loop added at the right endpoint)
  takes to reach $0$ when started at $j:=d^\eps(x,y)/\eps$, i.e., at most $j\,\bigl(3+2\,\lfloor
  \diam(X)/\eps \rfloor-j\bigr)$.  Hence
\begin{equation*}
  \begin{aligned}
&  \E\Bigl[F(x_\tau)+\sum_{i=0}^{\tau-1} \eps^2
    f(x_i)\Bigr]
\\ &\qquad
\ge  F(y) - \bigl(d^\eps(x,y) + 2\,\eps\bigr)\,\Lip_Y^\eps F  + j\,\bigl(3+2\,\lfloor
  \diam(X)/\eps \rfloor-j\bigr) \,\eps^2 \min(0,\inf f)
\\ &\qquad
  \ge  F(y) - \bigl(d^\eps(x,y) + 2\,\eps\bigr)\,\Lip_Y^\eps F  - d^\eps(x,y)\,\bigl(2\,\eps+2\,
  \diam(X)\bigr) \sup|f|.
\end{aligned}
\end{equation*}
This gives the desired lower bound on $v^\eps(x)$.
The symmetric argument gives the upper bound for $w^\eps(x)$.
\end{proof}

Next, we show that the lower bound $v^\eps$ on $u_\I^\eps$ is a good lower bound.

\begin{lemma}
  \label{vugap}
  Suppose $F$ is Lipschitz, and either
\begin{enumerate}
\item $f=0$ everywhere, or
\item $f$ is uniformly continuous, $X$ has finite diameter and\/ $\inf |f|>0$.
\end{enumerate}
  Then $\| u_\I^\eps - v^\eps\|_\infty\to 0$ as
  $\eps\to0$.  If $f$ is also Lipschitz, then $\| u_\I^\eps -
  v^\eps\|_\infty = O(\eps)$.
\end{lemma}

Note that since $X$ is assumed to be a length space,
the assumptions imply that $\sign(f)$ is constant and $\sup|f|<\infty$.

\begin{proof}
In order to prove that $u_\I^\eps$ is not much larger than $v^\eps$,
consider a strategy $S_\I$ for player~\I in ordinary $\eps$-tug-of-war,
which achieves an expected payoff of at least $u_\I^\eps-\eps$
against any strategy for player~\II.
We shortly describe a modified strategy
$S^\favored_\I$ for player~\I playing the \II-favored game, which
does almost as well as $S_\I$ does in the ordinary game.
To motivate $S^\favored_\I$, observe that a turn in
the \II-favored $\eps$-tug-of-war can alternatively be described as
follows. Suppose that the position at the end of the previous turn is
$x$. First, player~\I gets to make a move to an arbitrary point
$z$ satisfying $d(x,z)<\eps$.
Then a coin is tossed. If player~\II wins the toss,
she gets to make two steps from $z$, each of distance less than $\eps$.
Otherwise, \II gets to move to an arbitrary point in $B_{2\eps}(z)\cap Y$,
but only if the latter set is nonempty. This completes the turn.
The strategy $S^\favored_\I$ is based on the idea that after a win in the
coin toss by~\II, player~\I may use his move to reverse one of the
two steps executed by~\II (provided $Y$ has not been reached).

As strategy $S^\favored_\I$ is playing the \II-favored game against
player \II, it keeps track of a virtual ordinary game.
At the outset, the \II-favored game is in state
$x^\favored_0$, as is the virtual game.
As long as the virtual and the favored game have
not ended, each turn in the
favored game corresponds to a turn in the virtual game,
and the virtual game uses the same coin tosses as the favored
game.
In each such turn $t$, the target $z^\favored_t$ for player~\I in the favored game
is the current state $x_t$ in the ordinary game. If player~\I
wins the coin toss, then the new game
state $x^\favored_{t+1}$ in the favored game is his current target $z^\favored_t$,
which is the state of the virtual game $x_t$, and the new state of the
virtual game $x_{t+1}$ is chosen according to strategy $S_\I$
applied to the history of the virtual game.
If player~\II wins the toss and chooses the new state of the
favored game to be $x^\favored_{t+1}$, where
necessarily $d(x^\favored_{t+1},x_t)<2\,\eps$,
then in the virtual game the virtual player~\II chooses the new state
$x_{t+1}$ as some point satisfying $d(x_{t+1},x_t)<\eps$
and $d(x_{t+1},x^\favored_{t+1})<\eps$.
Induction shows that
$d(x^\favored_t,x_t)<\eps$ as long as both games are running,
and thus the described moves are all legal.

If at some time the virtual game has terminated, but the
favored game has not,
we let player~\I continue playing the favored game
by always pulling towards the final state of the virtual game.
If the favored game has terminated, for the sake of comparison,
we continue the virtual game, but this time let player~\II pull
towards the final state of the favored game and let
player~\I continue using strategy $S_\I$.

Let $\tau^\favored$ be the time at which the favored game has ended,
and let $\tau$ be the time at which the ordinary virtual game has ended.
By Lemma~\ref{vwgap}, if $\tau<\tau^\favored$, then
the conditioned expectation of the
remaining running payoffs and final payoff to player~\I in the
favored game after
time $\tau$, given what happened up to time $\tau$, is at least $F(x_{\tau})-O(\eps)$
(here the implicit constant may depend on $\diam X$, $\Lip_YF$ and $\sup|f|$).
Likewise, $u_I^\eps\le w^\eps$ from Lemma~\ref{epsilongameinequalities} and
the second inequality from Lemma~\ref{vwgap} show that
if $\tau^\favored<\tau$, then the conditioned expectation of the remaining
running and final payoffs in the virtual game is at most
$F(x^\favored_{\tau^\favored})+O(\eps)$.

Set $\lambda=\lambda_\eps:= \sup\bigl\{|f(x)-f(x')|:x,x'\in X\setminus Y,d(x,x')<2\,\eps\bigr\}$.
Then $\lambda = O(\eps)$ if $f$ is Lipschitz and $\lim_{\eps\to 0}\lambda_\eps =0$
if $f$ is uniformly continuous.
At each time  $t<\tau\wedge\tau^\favored$, since $z^\favored_t=x_t$,
the running payoff in the virtual game and in the favored game
differ by at most $\eps^2\,\lambda$.
Lemmas~\ref{epsilongameinequalities} and~\ref{vwgap} show
that there is a constant $C$, which may depend on $X,Y,f$ and $F$, but not on $\eps$,
such that
$-C\le v^\eps\le u_\I^\eps\le C$.
Thus, in the case where $\sup f<0$, since $S_I$ guarantees a payoff
of at least $u_\I^\eps(x_0)-\eps$, we have $\E[\tau^\favored\wedge\tau]=O(\eps^{-2})$.
Assume that player~\II plays the \II-favored game
(up to time $\tau\wedge\tau^\favored$) using a strategy $S^\favored_\II$
such that the expected payoff to player~\I who
uses $S^\favored_\I$ is at most $v_\eps+\eps$.
Then, in the case where $\inf f>0$,
 we will have $\E[\tau\wedge\tau^\favored]=O(\eps^{-2})$, again.
There is a strategy $S_\II$ for player~\II in the ordinary game,
which corresponds to the play of player~\II in the virtual game
when player~\I uses $S_I$ and $S^\favored_I$
and player~\II uses $S^\favored_\II$ in the favored game.
(The description of the virtual game defines $S_\II$ for some
game histories, and we may take an arbitrary extension of this
partial strategy to all possible game histories.)
The above shows that the expected payoff for player~\I
in the ordinary game when player~\I uses $S_\I$ and
player~\II uses $S_\II$ differs
from the expected payoff for player~\I
in the favored game when player~\I uses $S^\favored_\I$
and player~\II uses $S^\favored_\II$ by at most
$O(\eps)+\lambda\,\eps^2\,\E[\tau\wedge\tau^\favored]\le O(\eps+\lambda)$.
Thus $u^\eps_\I\le v^\eps+O(\eps+\lambda)$.
Since $u^\eps_\I\ge v^\eps$, the proof is now complete.
\end{proof}

Hence under the assumption of Lemma~\ref{vugap}, $\lim_{n\to\infty}
u_\I^{\eps 2^{-n}} = v^{\eps/2^\infty}$.

\subsection{Dyadic limits satisfy quadratic comparison}

We start by showing that $u_\I^\eps$ almost satisfies $(-2f)$-quadratic comparison from above.

\begin{lemma} \label{gcomparisonbound}
Let $\eps>0$,
  let $V$ be an open subset of $X \smallsetminus Y$ and write $V_\eps
  = \{x: \overline{B_\eps(x)} \subset V \}$.
  Suppose that $\phi(x)=Q(d(x,z))$
  is a quadratic distance function that is \increasing on $V$, where
  $Q(r)=ar^2+br+c$ satisfies
\begin{equation} \label{abound}
a \leq
-\sup_{x \in V_\eps}f(x).
\end{equation}
  Also suppose that
  $\sup_{V_\eps} f\ge 0$ or $\diam(V)<\infty$.
If the value function $u_{\I}^\eps$ for player~\I in
$\eps$-tug-of-war satisfies $u^\eps_{\I} \leq \phi$ on $V \backslash V_\eps$,
then $u^\eps_{\I} \leq \phi$ on $V_\eps$.
\end{lemma}

\begin{proof}
  Fix some $\delta>0$.
  Consider the strategy for player~\II that from a state $x_{k-1} \in
  V_\eps$ at distance $r=d(x_{k-1},z)$ from $z$ pulls to state $z$ (if
  $r<\eps$) or else moves to reduce the distance to $z$ by ``almost''
  $\eps$ units, enough to ensure that $Q(d(x_k,z)) < Q(r-\eps) +
  \delta 2^{-k}$.  If $r<\eps$, then $z \in V$ whence $Q(t)=at^2+c$
  with $a \ge 0$.  In this case, if \II wins the toss, then $x_k=z$,
  whence $\phi(x_k)=Q(0) \le Q(r-\eps)$. Thus for all $r \ge 0$,
  regardless of what strategy player~\I adopts,
  \begin{equation*}
  \E[\phi(x_k) \mid x_{k-1}] - \delta 2^{-k-1} \le \frac{Q(r+\eps)+Q(r-\eps)}{2}
  =Q(r)+a\,\eps^2
  = \phi(x_{k-1})+a\,\eps^2\,.
  \end{equation*}
  Setting
  $\tau:=\inf \{k : x_k \notin V_\eps\}$, we conclude that
  $M_k:=\phi(x_{k\wedge\tau})-a\,\eps^2(k\wedge\tau)+\delta\,2^{-k}$ is a supermartingale.

  Suppose that player~\I uses a strategy with expected payoff larger than $-\infty$.
  (If there is no such strategy, the assertion of the Lemma is obvious.)
  Then $\tau<\infty$ a.s.  We claim that
  \begin{equation}\label{supmart}
  \E[M_\tau] \le M_0 \,.
  \end{equation}
  Clearly this holds if $\tau$ is replaced by $\tau \wedge k$. To pass
  to the limit as $k \to \infty$, consider two cases.
  \begin{itemize}
  \item If $a\le 0$, then since $\phi$ is \increasing on $V$, it is also bounded from below on $V$.
    Consequently, $M_{k}$ is a supermartingale bounded below, so \eqref{supmart} holds.
  \item
    If $a>0$, then $\sup_{V_\eps} f<0$, by~\eqref{abound}.
    By assumption therefore $\diam V<\infty$, which implies $\sup_V|\phi|<\infty$.
    If $\E[\tau]=\infty$, we get $\E[M_\tau]=-\infty$, and hence~\eqref{supmart} holds.
    On the other hand, if $\E[\tau]<\infty$, then dominated convergence gives~\eqref{supmart}.
  \end{itemize}
  Since $u^\eps_\I(x_\tau) \le \phi(x_\tau)$, we deduce
  that
  $$ u^\eps_\I(x_0) \le \sup_{S_\I} \E\Bigl[\phi(x_\tau)+\sum_{t=0}^{\tau-1} f(x_t)\Bigr]
  \overset{\eqref{abound}}{\le} \sup_{S_\I}\E[ M_\tau] \le M_0 = \phi(x_0)+\delta\,.
  $$
  where $S_\I$ runs over all possible strategies for player~\I
  with expected payoff larger than $-\infty$. Since $\delta>0$ was arbitrary,
  the proof is now complete.
\end{proof}

In order for $u_\I^\eps$ to satisfy $(-2f)$-quadratic comparison (from above), we
would like to know that if $u_\I^\eps \leq \phi$ on the boundary of an
open set, then this (almost) holds in a neighborhood of the boundary,
so that we can apply the above lemma.  To do this we prove a uniform
Lipschitz lemma:

\begin{lemma}[Uniform Lipschitz]
  \label{uniformlylipschitz}
Suppose that $F$ is Lipschitz, and
either
      \begin{enumerate} \item $f=0$ everywhere,\, or  \item $|f|$ is bounded above and $X$ has finite diameter.
      \end{enumerate}
Then for each $\eps\in(0,\diam X)$, $u_\I^\eps$ and $u_\II^\eps$ are
Lipschitz on $X$ w.r.t.\ the metric $d^\eps$ with Lipschitz constant
depending only on $\diam X,\sup|f|$ and $\Lip_Y F $.
\end{lemma}
\begin{proof}
By symmetry, it suffices to prove this for $u_\I^\eps$.
Set $L:= 3\,\Lip^\eps_Y F + 4\, \diam X\sup|f|$.
Let $x,y\in X$ be distinct. If $x,y\in Y$, then
$\bigl|u_\I^\eps(x)-u_\I^\eps(y)\bigr|=\bigl|F(x)-F(y)\bigr|\le \Lip^\eps_Y F \,d^\eps(x,y)$.
If $x\in X\setminus Y$ and $y\in Y$,
Lemmas~\ref{vwgap} and~\ref{epsilongameinequalities} give
\begin{equation}\label{xinyout}
\bigl|u_\I^\eps(x)-u_\I^\eps(y)\bigr|=
\bigl|u_\I^\eps(x)-F(y)\bigr|\le L\,d^\eps(x,y)\,,\qquad x\in X\setminus Y,\quad y\in Y\,.
\end{equation}
Now suppose $x,y\in X\setminus Y$.
Set $Y^*=Y\cup\{y\}$, $F^*=F$ on $Y$ and $F^*(y)=u_\I^\eps(y)$.
Then, clearly, the value of $u_\I^\eps(x)$
for the game where $Y$ is replaced by $Y^*$ and $F$ is
replaced by $F^*$ is the same as for the original game.
By~\eqref{xinyout}, we have $\Lip^\eps_{Y^*}(F^*)\le L$.
Consequently,~\eqref{xinyout} gives
$$
\bigl|u_\I^\eps(x)-u_\I^\eps(y)\bigr|
=
\bigl|u_\I^\eps(x)-F^*(y)\bigr| \le
\bigl(3\,L + 2\, (\eps +\diam X)\,\sup|f|\bigr)\,
d^\eps(x,y)
\le
4\,L\,
d^\eps(x,y)\,,
$$
which completes the proof.
\end{proof}

\begin{lemma}
  \label{dyadic-comparison}
  Suppose $F$ is Lipschitz and $\inf F>-\infty$, and either (1)
  $f=0$ identically or else (2) $\inf f>0$, $f$ is uniformly
  continuous, and $X$ has finite diameter.  Then the subsequential
  limit $v^{\eps/2^\infty}=\lim_{n\to\infty} v^{\eps 2^{-n}}$ satisfies
  $(-2f)$-quadratic comparison on $X\smallsetminus Y$.
\end{lemma}

\begin{proof}
  By Theorems~\ref{uniqueinfharmongraph} and
  \ref{uniquesolutionongraph}, $u_\I^\eps=u_\II^\eps$, so from
  Lemma~\ref{vugap} we have $\|w^\eps-v^\eps\|_\infty\to 0$ as
  $\eps\to 0$.  Thus $\lim_{n\to\infty} w^{\eps 2^{-n}} = v^{\eps/2^\infty}$.

  Note that the hypotheses imply that $|f|$ is bounded.
  Consider an open $V \subset X\smallsetminus Y$ and an
  \increasing quadratic distance function $\phi$ on $V$ with quadratic
  term $a \leq -\sup_{y \in V} f(y) $, such that $\phi \geq
  v^{\eps/2^\infty}$ on $\partial V$.  We must show that $\phi \ge
  v^{\eps/2^\infty}$ on $V$.  Since $\|w^\eps-v^\eps\|_\infty\to 0$, we have
  by Lemma~\ref{epsilongameinequalities}
  that $v^{\eps 2^{-n}}$ converges uniformly to $v^{\eps/2^\infty}$.  So for
  any $\delta >0$, if $n\in\N$ is large enough, then $v^{\eps 2^{-n}}
  \le \phi+\delta$ on $\partial V$.
  Note also that $\phi$ is necessarily uniformly continuous on $\overline V$.
  (If $\diam V<\infty$ or $a=0$, this is clear. Otherwise, $f=0$ and $a<0$.
  However, $a<0$ implies that $\diam V<\infty$, since $\phi$ is \increasing on $V$.)
  Hence, by the uniform Lipschitz
  lemma~\ref{uniformlylipschitz},
  $u_\I^{\eps 2^{-n}} \le \phi+2\,\delta$ on $V
  \smallsetminus V_{\eps 2^{-n}} $ for all sufficiently large $n\in \N$,
  where we use the notations of Lemma~\ref{gcomparisonbound}.
  By that lemma $u_\I^{\eps 2^{-n}} \le \phi+2\,\delta$ on all of $V$.
  Letting $n\to\infty$ and $\delta\to 0$ shows that $v^{\eps/2^\infty}$
  satisfies $(-2f)$-quadratic comparison from above.

  To prove quadratic comparison from below, note that the only assumptions
  which are not symmetric under exchanging the roles of the players are
  $\inf F>-\infty$ and $\inf f>0$.
  However, we only used these assumptions to  prove
  $w^{\eps/ 2^\infty} = v^{\eps/2^\infty}$. Consequently, comparison from below
  follows by symmetry.
\end{proof}

\subsection{Convergence}

\begin{lemma}
  \label{viscosity-supermartingale}
  Suppose that $v$ is continuous and satisfies $(-2f)$-quadratic comparison from below
  on $X\smallsetminus Y$,
  and that $f$ is locally bounded below. Let $\delta>0$.  Then in \II-favored
  $\eps$-tug-of-war, when player~\I (using any strategy)
  targets point $z_i$ on step $i$,
  player~\II may play to make $M_{t\wedge\tau_\eps}$ a supermartingale, where
  $$ M_t := v(x_t)+\eps^2\sum_{i=1}^{t} \inf_{y\in B_{2\eps}(z_i)} f(y) +
  \delta 2^{-t}$$
  and $\tau_\eps:=\inf\{t:d(x_t,Y)<3\,\eps\}$.
\end{lemma}

\begin{proof}
  Let $z=z_t$ be the point that player~\I has targeted at time $t$.
  Assume that $t<\tau_\eps$.
  We define the following:
\begin{enumerate}
\item $\alpha := \inf \{ f(x): x \in B_{2\eps}(z)\}$,
\item $A := \inf \{v(x): x \in B_{2\eps}(z)\}$ (the
infimum value of $v$ that \II can guarantee if \II wins the coin toss),
\item $\beta := \frac{v(z) + A}{2} + \alpha\,\eps^2$, and
\item $Q(r) := -\alpha r^2 + \frac{A-v(z)+4\alpha\eps^2}{2\eps} r + v(z)$\ \  ($Q(0)=v(z)$, $Q(2\eps)=A$, $Q(\eps)=\beta$, and $Q''=-2\alpha$).
\end{enumerate}
Player~\II can play so that $\E\bigl[v(x_t)\bigm|\text{$z_t$ and all prior
  events}\bigr] \leq (v(z)+A)/2+ \delta 2^{-t}$, i.e., so that
$\E\bigl[M_t\bigm| z_t \text{ and prior events}\bigr]-M_{t-1} \leq \beta-v(x_{t-1})$.  We will show that
whenever $d(x,z)<\eps$ we have $v(x) \geq \beta$, and then it will
follow that $M$ is a supermartingale.  There are two cases to check,
depending on whether $\alpha>0$ or $\alpha\leq 0$:

Suppose $\alpha \leq 0$. Note that  $A\leq v(z)$.
If $v(z)=A$, the inequality $v(x)\ge\beta$ on $B_\eps(z)$ follows
from $\alpha\le 0$ and the definition of $\beta$.
Assume therefore that $v(z)>A$. Then
$Q'(\eps)=(A-v(z))/(2\eps)< 0$.  Thus $Q$ is decreasing on
$[0,\eps]$.  Let $r_0\in [\eps,2\,\eps]$ be the point where
$Q$ attains its minimum in $[\eps,2\,\eps]$.
Then $Q$ is decreasing on $[0,r_0]$.
Set $V:=B_{r_0}(z)\smallsetminus \{z\}$.
We have $v(z)=Q(0)=Q\bigl(d(z,z)\bigr)$
and for $x\in\partial B_{r_0}(z)$ we have
$v(x)\geq A=Q(2\,\eps) \geq Q(r_0)=Q\bigl(d(x,z)\bigr)$.
Thus, $v(x)\ge Q\bigl(d(x,z)\bigr)$ for $x\in\partial V$.
Since $v$ satisfies $(-2f)$-quadratic comparison from below,
and $Q''=-2\,\alpha\geq \sup_{x\in V} -2\,f(x)$,
we get $v(x)\geq Q\bigl(d(x,z)\bigr)$ for $x\in V$.
In particular, for $x\in B_\eps(z)$ one has
$v(x)\ge Q\bigl(d(x,z)\bigr)\ge Q(\eps)=\beta$.

Now suppose $\alpha>0$.
The function $L_0(x) = -\alpha\, d(x,z)^2 + \alpha(2\,
\eps)^2 + A$ is a lower bound for $v$ on $\partial B_{2\eps}(z)$, and hence
applying $(-2f)$-quadratic comparison in $B_{2\eps}(z)$ with $L_0(\cdot)$
gives $v(z) \geq L_0(z) = A + 4\,\alpha\,\eps^2$.
Therefore, $Q'(0)= (A-v(z) +
4\alpha\eps^2)/(2\eps)\leq 0$, which together with $Q''< 0$ implies
that $Q$ is decreasing on $[0,2\eps]$.  By applying $(-2f)$-quadratic comparison
on $B_{2\eps}(z)\setminus\{z\}$ we see that for $x\in
B_\eps(z)$ we have $v(x)\geq Q(d(x,z))\geq Q(\eps)=\beta$.
\end{proof}

\begin{lemma}
  \label{vepsilonlowerboundsviscosity}
Suppose that $F$ is Lipschitz, and either
\begin{enumerate} \item $f=0$ everywhere,\, or  \item $|f|$ is bounded above and $X$ has finite diameter.
\end{enumerate}
  Also suppose that $v:X\to\R$ is continuous,
  satisfies $(-2f)$-quadratic comparison from below on $X\smallsetminus Y$,
  $v\ge F$ on $Y$, and\/ $\inf v>-\infty$.  Then $v^\eps \leq v$ for
  all $\eps$.
\end{lemma}

\begin{proof}
  The idea is for player~\II to make the $M$ defined in
  Lemma~\ref{viscosity-supermartingale} a supermartingale, but we need
  to pick a stopping time $\tau$ such that $\E[M_\tau]\leq M_0$ while
  $v^\eps(x_\tau)$ is unlikely to be much larger than $v(x_\tau)$.
  Let $W:=\{x\in X:d(x,Y)\ge 3\,\eps\}$ and
  let $\tau_\eps$ be defined as in Lemma~\ref{viscosity-supermartingale};
  that is $\tau_\eps:=\inf\{t\in\N: x_t\notin W\}$.
  Set $\lambda_\eps:=\sup_{X\setminus W}(v^\eps-v)$,
  and let $\delta>0$.
  We first show
\begin{equation}\label{e.firstshow}
v^\eps\le v+\lambda_\eps+\delta\,.
\end{equation}

  In the case that $f=0$, the supermartingale $M$ is bounded below, so we
  can choose $\tau=\tau_\eps$. Player~\I is compelled to ensure $\tau<\infty$.
  Conditional on the game up to time $\tau$,
  player~\I cannot guarantee a conditional expected payoff better than
  $v^\eps(x_\tau)+\delta\,2^{-\tau}\le M_\tau+\lambda_\eps$.
  Since $\E[M_\tau]\leq M_0 = v(x_0)+\delta$, we get~\eqref{e.firstshow},
   as desired.

  In the case $f\neq 0$, we let $\tau_n:=\tau_\eps\wedge n$, where $n\in\N$.
  Then
  $\E[M_{\tau_n}] \leq M_0$.
Note that $\sup v^\eps<\infty$, follows from $\diam X<\infty$, $v^\eps\le u_\I^\eps$
and Lemma~\ref{uniformlylipschitz}.
Suppose player~\II makes $M$ a
  supermartingale up until time $\tau_n$.  Given
play until time $\tau_n$, player~\II may make
sure that the conditional expected payoff to player~\I is at most
  $$\delta\,2^{-\tau_n}+v^\eps(x_{\tau_n})+\eps^2\sum_{i=1}^{\tau_n} \inf_{y\in
      B_{2\eps}(z_i)} f(y)\,.
$$
Taking expectation and separating into cases in which $x_{\tau_n}\in W$ or not, we get
  $$v^\eps(x_0)\le \E\Bigl[v(x_{\tau_n})+\eps^2\sum_{i=1}^{\tau_n} \inf_{y\in
      B_{2\eps}(z_i)} f(y)\Bigr] + \lambda_\eps+\delta\,2^{-\tau_n}+
  \Pr[x_{\tau_n}\in W]\,\sup_W \bigl(v^\eps(x)-v(x)\bigr)\,.
$$
Since $
  \E[M_{\tau_n}]\leq M_0 = v(x_0)+\delta$, this gives
  $$v^\eps(x_0)\le
      v(x_0)+\delta+ \lambda_\eps+ \Pr[x_{\tau_n}\in W]\,\sup_W \bigl(v^\eps(x)-v(x)\bigr)\,.
$$
  The first term above is $<
  \E[M_{\tau_n}]\leq M_0 = v(x_0)+\delta$, independent of $n$.
  Player~\I is compelled to play a strategy that ensures $\tau_\infty$
  is finite a.s., since otherwise the payoff is $-\infty < v(x)$.
  With such a strategy, $\Pr[x_{\tau_n}\in W]\to 0$ as
  $n\to\infty$, and since $v$ is bounded below and $\sup v^\eps<\infty$, the last summand
  tends to $0$ as $n\to\infty$. Thus, we get~\eqref{e.firstshow} in this case as well.

Next, we show that $\limsup_{\eps\searrow 0}\lambda_\eps\le 0$.
Let $y\in Y$, and set $Q(r)=a\,r^2+b\,r+c$, where $a:=\sup|f|$,
$b<b^*:= - 2\,\sup|f|\,\diam X-\Lip F$,
and $c:= F(y)$.  Let $\phi(x):=Q\bigl(d(x,y)\bigr)$.
Then $\phi(y')\le F(y')\le v(y')$ for $y'\in Y$.
Since $v$ satisfies comparison from below and $\phi$ is \decreasing,
we get $v(x)\ge \phi(x)$ on $X$.
Thus, if $x\in \overline {B_{3\eps}(y)}$,
then $v(x)\ge F(y) +b^*\,d(x,y)\ge F(y)+3\,b^*\,\eps$.
In conjunction with Lemma~\ref{uniformlylipschitz}, this implies
that $\limsup_{\eps\searrow 0}\lambda_\eps\le 0$.
Choosing $\delta=\eps$ and taking
$\eps$ to $0$ therefore gives in~\eqref{e.firstshow} $\limsup_{\eps\searrow0}v^\eps\le v$.
However, the inequality $v^{2\eps}\le v^\eps$ from Lemma~\ref{epsilongameinequalities}
implies $v^\eps\le \limsup_{\eps'\searrow0}v^{\eps'}\le v$,
completing the proof.
\end{proof}

\begin{proof}[Proof of Theorem~\ref{continuumvalueexists}]
First, suppose that $F$ is Lipschitz.
Let $\eps,\eps'>0$.
Let $v:=\lim_{n\to\infty} v^{\eps2^{-n}}$
and $v':= \lim_{n\to\infty} v^{\eps'2^{-n}}$.
We know that these limits exist from Lemma~\ref{epsilongameinequalities}.
Lemma \ref{vugap} tells us that
$\|u_I^\eps-v^\eps\|_\infty\to 0$ as $\eps\to0$.
Since the assumptions of that Lemma are player-symmetric,
we likewise get $\|u_\II^\eps-w^\eps\|_\infty\to 0$.
{}From Theorems~\ref{uniqueinfharmongraph} and~\ref{uniquesolutionongraph}
(possibly with the roles of the players interchanged)
we know that $u_\I^\eps=u_\II^\eps$, and hence the above gives
$\|v^\eps-w^\eps\|_\infty\to 0$.
By Lemma~\ref{epsilongameinequalities}, $v^\eps\le v\le w^\eps$
and  $v^\eps\le u_\I^\eps\le w^\eps$,
and so we conclude that $\|v^\eps-v\|_\infty\to 0$
and $\|u_\I^\eps- v\|_\infty\to 0$.
Note that the assumptions imply that $\sup|f|<\infty$.
Therefore, from Lemma~\ref{uniformlylipschitz} and
$\|u_\I^\eps-v\|_\infty\to 0$ we conclude that $v$ is
Lipschitz.
Precisely the same argument gives
$\|v^\eps-v\|_\infty=O(\eps)$ if $f$ is also assumed to be Lipschitz,
and
similar estimates also hold for $\|v^{\eps'}-v'\|_\infty$.
Note that the assumptions imply that $F$ is bounded if $f\ne 0$.
Lemma~\ref{dyadic-comparison}
(applied possibly with the roles of the players reversed)
tells us that $v'$ satisfies
$(-2f)$-quadratic comparison on $X\setminus Y$.
Clearly $\inf v'>-\infty$.
(If $f=0$ identically, then $\inf v'\ge\inf F$,
while if $\diam X<\infty$, we may use the fact that $v'$ is Lipschitz.)
Thus, Lemma~\ref{vepsilonlowerboundsviscosity}
implies
that $v^\eps\le v'$. Consequently, $v\le v'$.
By symmetry $v'\le v$, and hence $v=v'$.
This completes the proof in the case where $F$ is Lipschitz.

Using the result of Lemma~\ref{extensionconditions} below,
we know that for every $\delta>0$ there is a Lipschitz $F_\delta:Y\to\R$
such that $\|F-F_\delta\|_\infty<\delta$.
Then  the Lipschitz case applies to the functions $F_\delta\pm \delta$ in place
of $F$. Since the game value $u^\eps$ for $F$ is bounded between
the corresponding value with $F_\delta+\delta$ and $F_\delta-\delta$,
and the latter two values differ by $2\delta$,
the result easily follows.
\end{proof}

\begin{lemma}\label{extensionconditions}
Let $X$ be a length space, and let $F:Y\to\R$ be defined
on a nonempty subset $Y\subset X$.
The following conditions are equivalent:
\begin{enumerate}
\item
$F$ is uniformly continuous and $\sup \{ (F(y)-F(y'))/\max\{1,d(y,y')\}:y,y' \in Y\} < \infty$.
\item $F$ extends to a uniformly continuous function on $X$.
\item There is a sequence of Lipschitz functions on $Y$ tending to $F$ in
  $\|\cdot\|_\infty$.
\end{enumerate}
\end{lemma}
\begin{proof}
We start by assuming 1 and proving 2.
Let $ \phi(\delta):= \sup\bigl\{F(y)-F(y'):d(y,y')\le\delta,\,y,y'\in Y\bigr\}$,
and $\tilde\phi(t):=\sup\bigl\{t\,\phi(\delta)/\delta: \delta\ge t\bigr\}$.
We now show that $\lim_{t\searrow0}\tilde\phi(t)=0$. Let $\eps>0$,
and let $\delta_\eps>0$ satisfy $\phi(\delta_\eps)<\eps$. Such a $\delta_\eps$
exists because $F$ is uniformly continuous.
Condition 1 implies that $M:=\sup\{\phi(\delta)/\delta:\delta\ge\delta_\eps\}<\infty$.
For $t<(\eps/M)\wedge \delta_\eps$, we have
$$
\tilde\phi(t)=
\sup\bigl\{t\,\phi(\delta)/\delta:\delta_\eps\ge \delta\ge t\bigr\}
\vee
\sup\bigl\{t\,\phi(\delta)/\delta:\delta> \delta_\eps\bigr\}
\le \phi(\delta_\eps)\vee(t\,M)\le \eps\,,
$$
which proves that $\lim_{t\searrow0}\tilde\phi(t)=0$.
Two other immediate properties of $\tilde \phi$ which we will use
are that
$F(y)-F(y')\le\phi\bigl(d(y,y')\bigr)\le \tilde\phi\bigl(d(y,y')\bigr)$
holds for $y,y'\in Y$ and
$\tilde\phi(s\,t)\ge s\,\tilde\phi (t)$ when $s\in[0,1]$ and $t\ge0$.
It follows that $\tilde\phi$ is subadditive: $\tilde\phi(a)+\tilde\phi(b)\ge
a\,\tilde\phi(a+b)/(a+b)+b\,\tilde\phi(a+b)/(a+b)=\tilde\phi(a+b)$
for $a,b\ge 0$.

Now for $x\in X$ set $u(x):=\inf\bigl\{F(y)+\tilde\phi(d(y,x)):y\in Y\bigr\}$.
Then $u=F$ on $Y$. We now prove
\begin{equation}\label{e.uu}
u(x)-u(x')\le \tilde\phi\bigl(d(x,x')\bigr)
\end{equation}
for $x,x'\in X$.
Indeed, let $y'\in Y$.
Then
\begin{multline*}
u(x)
- F(y') - \tilde \phi\bigl(d(y',x')\bigr)
\le
F(y')+\tilde\phi\bigl(d(y',x)\bigr)
- F(y') - \tilde \phi\bigl(d(y',x')\bigr)
\\
=
\tilde\phi\bigl(d(y',x)\bigr)
 - \tilde \phi\bigl(d(y',x')\bigr).
\end{multline*}
Consequently, subadditivity and monotonicity of $\tilde\phi$ gives
$$
u(x)
- F(y') - \tilde \phi\bigl(d(y',x')\bigr)
\le \tilde\phi\bigl(\bigl|d(x,y')-d(x',y')\bigr|\bigr)
\le \tilde\phi\bigl(d(x,x')\bigr).
$$
Taking the supremum over all $y'\in Y$ then implies~\eqref{e.uu}.
Therefore, $u$ is uniformly continuous and 2 holds.

\smallskip
We now assume 2 and prove 3. Let $u:X\to\R$ be a uniformly continuous
extension of $F$ to $X$. It clearly suffices to approximate $u$
by Lipschitz functions on $X$ in $\|\cdot\|_\infty$.
For $x\in X$ let $u_L(x):= \inf\bigl\{u(x')+L\,d(x,x'):x'\in X\bigr\}$.
The same argument which was used above to prove~\eqref{e.uu} now shows
that $\Lip(u_L)\le L$. Clearly, $u_L(x)\le u(x)$.
Let $\phi(t):=\sup\bigl\{u(x)-u(x'):d(x,x')\le t\bigr\}$ for $t\ge 0$.
Since $X$ is a length space, $\phi$ is subadditive.
Let $t>0$ and $k:=\lfloor d(x,x')/t\rfloor$.
Then
$$
u(x)-u(x')-L\,d(x,x')\le
\phi\bigl(d(x,x')\bigr)-L\,d(x,x')
\le
\phi\bigl((k+1)\,t\bigr)-L\,k\,t\,.
$$
Taking the supremum over all $x'$ and using the subadditivity of $\phi$
therefore gives
$$
u(x)-u_L(x) \le \sup_{k\in\N}
\Bigl(
\phi\bigl((k+1)\,t\bigr)-L\,k\,t\Bigr)
\le \phi(t)+ \sup_{k\in\N}k\,\bigl( \phi(t)-L\,t\bigr) .
$$
Therefore, $u(x)-u_L(x)\le \phi(t)$ once $L>\phi(t)/t$.
Since $\inf_{t>0} \phi(t)=0$ and $u_L\le u(x)$, this proves 3.

\smallskip
The passage from 3 to 1 is standard, and therefore omitted.
This concludes the proof.
\end{proof}

\begin{proof}[Proof of Theorem~\ref{deltainftyuisf}]
Lemma~\ref{extensionconditions} tells us that $F$ extends to a uniformly
continuous function on $X$ and that it can be approximated in $\|\cdot\|_\infty$
by Lipschitz functions.
We know from Theorem~\ref{continuumvalueexists} that the continuum value
$u$ exists and is uniformly continuous.
Suppose first that $F$ is Lipschitz.
Lemma~\ref{dyadic-comparison} (applied possibly with the
roles of the players reversed) says that $u$ satisfies
$(-2f)$-quadratic comparison on $X\setminus Y$.
If $F$ is not Lipschitz, we may deduce the same
result by approximating $F$ from below by Lipschitz functions
and observing that a monotone nondecreasing
$\|\cdot\|_\infty$-limit of functions satisfying
$(-2f)$-quadratic comparison from above also
satisfies $(-2f)$-quadratic comparison from above,
and making the symmetric argument for comparison
from below.
Now suppose that $\tilde u$ is as in the second part of
the theorem.
This clearly implies $\inf\tilde u >-\infty$.
Now Lemma~\ref{vepsilonlowerboundsviscosity} gives
$u\le\tilde u$ (again, we may need to first approximate
$F$ by Lipschitz functions).
The uniqueness follows directly.
\end{proof}

\begin{proof}[Proof of Theorem~\ref{uniqueAM}]
If $x\in X\setminus Y$ and $y\in Y$,
then from Lemma~\ref{vwgap} it follows that
$\bigl|u(x)-u(y)\bigr|\le \Lip_YF\, d(x,y)$.
Let $U\subset X\setminus Y$ be open
and define $F^*(x)=u(x)$ for $x\in Y\cup\partial U$.
It is clear that the continuum value $u^*$ of
$(X,d,Y\cup\partial U,F^*,0)$ is the same as $u$,
since any player may first play a strategy
that is appropriate for the
$(X,d,Y\cup\partial U,F^*,0)$ game and once
$Y\cup\partial U$ is hit start playing a strategy
that is appropriate for the original game.
The above argument shows that
$\bigl|u(x)-u(y)\bigr|\le d(x,y)\,\Lip_{\partial U}u$
if $y\in \partial U$ and $x\in U$.
In particular, $\Lip_{\partial U\cup\{x\}}u=\Lip_{\partial U}u$,
which implies $\Lip_U u=\Lip_{\partial U}u$; that is, $u$ is
AM in $X\setminus Y$.

Now suppose that $u^*:X\to\R$ is an AM extension of $F$ and $\sup|F|<\infty$.
Lemma~\ref{AMiffCDF} tells us that $u^*$ satisfies comparison with
distance functions.
Observe that the proof of Lemma~\ref{vepsilonlowerboundsviscosity}
shows that in the case $f=0$ we may replace the hypothesis
that $v$ satisfies $0$-quadratic comparison from below
on $X\setminus Y$ by the hypothesis that $v$ satisfies
comparison with distance functions from below (since
only comparisons with distance functions are used in this case).
Thus, $u\le u^*$. Similarly, $u\ge u^*$, which implies
the required uniqueness statement and completes the proof.
\end{proof}

\section{Harmonic measure for $\Delta_\infty$} \label{section:porous}

Here, we present a few estimates of the $\infty$-harmonic measure $\omega_\infty$.
Before proving Theorem~\ref{t.hmeasure}, we consider the $\infty$-harmonic
of porous sets.
Recall that a set $S$ in a metric space $Z$ is $\alpha$-\textbf{porous}
if for every $r\in(0,\diam Z)$
every ball of radius $r$ contains
a ball of radius $\alpha\,r$ that is disjoint from $S$.
An example of a porous set is the ternary Cantor set in $[0,1]$.
We start with a general lemma in the setting of length spaces.

\begin{lemma}\label{game-plan}
  Suppose $X$ is a length space and $0\leq F\leq 1$ on the
  terminal states $Y$, where $F:Y\to\R$ is continuous.
  Let $S:=\supp F$.
  Suppose that for some integer $k$ and positive constants $\eps_0,d_{\min}$
  and $\gamma\in(0,1)$, for every $x\in X$ with $d(x,S)\geq
  d_{\min}$, there is a sequence of points $x=z_0,\ldots,z_k$ such that
  $z_k\in Y$, $F(z_k)=0$, and
 $$d(z_i,S) \geq 2\, d(z_i,z_{i-1}) + 2\,\eps_0 + \gamma\, d(z_0,S)$$
 for $i=1,\ldots,k$.  Then the AM extension $u$ of $F$
 satisfies $$0\leq u(x) \leq (1-2^{-k})^{\log_\gamma (d_{\min}/d(x,S))}.$$
\end{lemma}

\begin{proof}
  We will obtain bounds on $u^\eps(x)$ that are independent of $\eps$
  (as long as $\eps\in(0,\eps_0)$),
  and these yield bounds on $u(x)$.  Since $u^\eps\geq 0$, we need
  only give a good strategy for player~II to obtain an upper bound on
  $u^\eps(x)$.  The idea is for player~II to always have a ``plan''
  for reaching a terminal state at which $F$ is $0$ while staying away
  from the terminal states at which $F$ is nonzero.  A plan consists
  of a sequence of points $z_0,z_1,\ldots,z_k$, where $z_0$ is the
  game state when the plan was formed, $z_k\in Y$, and $F(z_k)=0$.  As
  soon as the game state reaches $z_i$, player~II starts pulling
  towards $z_{i+1}$.  If player~I is lucky and gets many moves,
  player~II may have to give up on the plan and form a new plan.  When
  tugging towards $z_i$, we suppose that player~II gives up and forms
  a new plan as soon as $d^\eps(x_t,z_i)=2 d^\eps(z_i,z_{i-1})$, and
  otherwise plays to ensure that $d^\eps(x_t,z_i)$ is a
  supermartingale.  While tugging towards $z_i$, the plan will be
  aborted at that stage with probability at most $1/2$, so with probability
  at least $2^{-k}$ the plan is never aborted and succeeds in reaching $z_k$.

  Suppose player~II aborts the plan at time $t$ while tugging towards
$z_i$, and forms
  a new plan starting at $z'_0 = x_t$.  Then $d(x_t,z_i)\leq d^\eps(x_t,z_i) = 2\, d^\eps(z_i,z_{i-1}) \leq 2\, d(z_i,z_{i-1}) + 2\,\eps$, so
$$
  d(x_t,S) \geq d(z_i,S) - d(x_t,z_i) \geq 
  d(z_i,S) - [2 d(z_i,z_{i-1})+2\eps] \geq \gamma\, d(z_0,S),
$$
  so the game state remains far from $S$.  Since player~II can
  always find a short plan (length at most $k$) when the distance from
  $S$ is at least $d_{\min}$, player~I gets to $S$ with probability
  at most $(1-2^{-k})^{\lceil\log_\gamma( d_{\min}/d(x_0,S))\rceil}$, which yields
  the desired upper bound.
\end{proof}

The reader may wish to check that the lemma  implies
$\lim_{\delta\searrow 0}\omega_\infty(A_\delta)=0$ in the setting of
Theorem~\ref{t.hmeasure} (i.e., when $X$ is the unit ball in
$\R^n$, $n>1$, $Y=\partial X$ and $A_\delta$ is a spherical cap of
radius $\delta$.)

We are now ready to state and prove an upper bound on the $\infty$-harmonic
measure of neighborhoods of porous sets.

\begin{theorem}
Let $X\subset \R^n$, $n>1$, be the closed unit ball and let $Y$ be its boundary,
the unit sphere.
Let $\alpha\in (0,1/2)$ and let $\delta>0$.
Let $S$ be an $\alpha$-porous subset of $Y$,
and let $S_\delta$ be the closure of the $\delta$-neighborhood
of $S$.
Then
$$
\omega^{(0,Y,X)}_\infty(S_\delta) \le  \delta^{\alpha^{O(1)}}.
$$
\end{theorem}

Of course, in the above, $S$ is $\alpha$-porous as a subset of $Y$.
(Every subset of $Y$ is $(1/3)$-porous as a subset of $X$.)

\begin{proof}
The plan is to use the lemma, of course.
Let $d_{\min}:= 2\,\delta/\alpha$.
Let $z_0\in X\setminus Y$, and suppose that $d_0:=d(z_0,S)\ge d_{\min}$.
Let $y_0\in Y$ be a closest point to $z_0$ on $Y$.
Inside $B_{d_0}(y_0)\cap Y$ there is a point
$y_1$ such that $B_{\alpha d_0}(y_1)\cap S=\emptyset$.
(See Figure~\ref{cantor}.)
Therefore, $d(y_1,S_\delta)\ge \alpha\, d_0-\delta \ge \alpha\,d_0/2$.
We define the sequence $z_j$ inductively, as follows.
If $d(z_j,S_\delta)\ge 3 \,d(z_j,Y)$, then we take $z_{j+1}$
to be any closest point to $z_j$ on $Y$.
Otherwise, let $z_{j+1}$ be the point
on the line segment from $z_j$ to $y_1$
whose distance from $z_j$ is $d(z_j,y_1)/10$.
\begin{figure}[h]
\centerline{\includegraphics{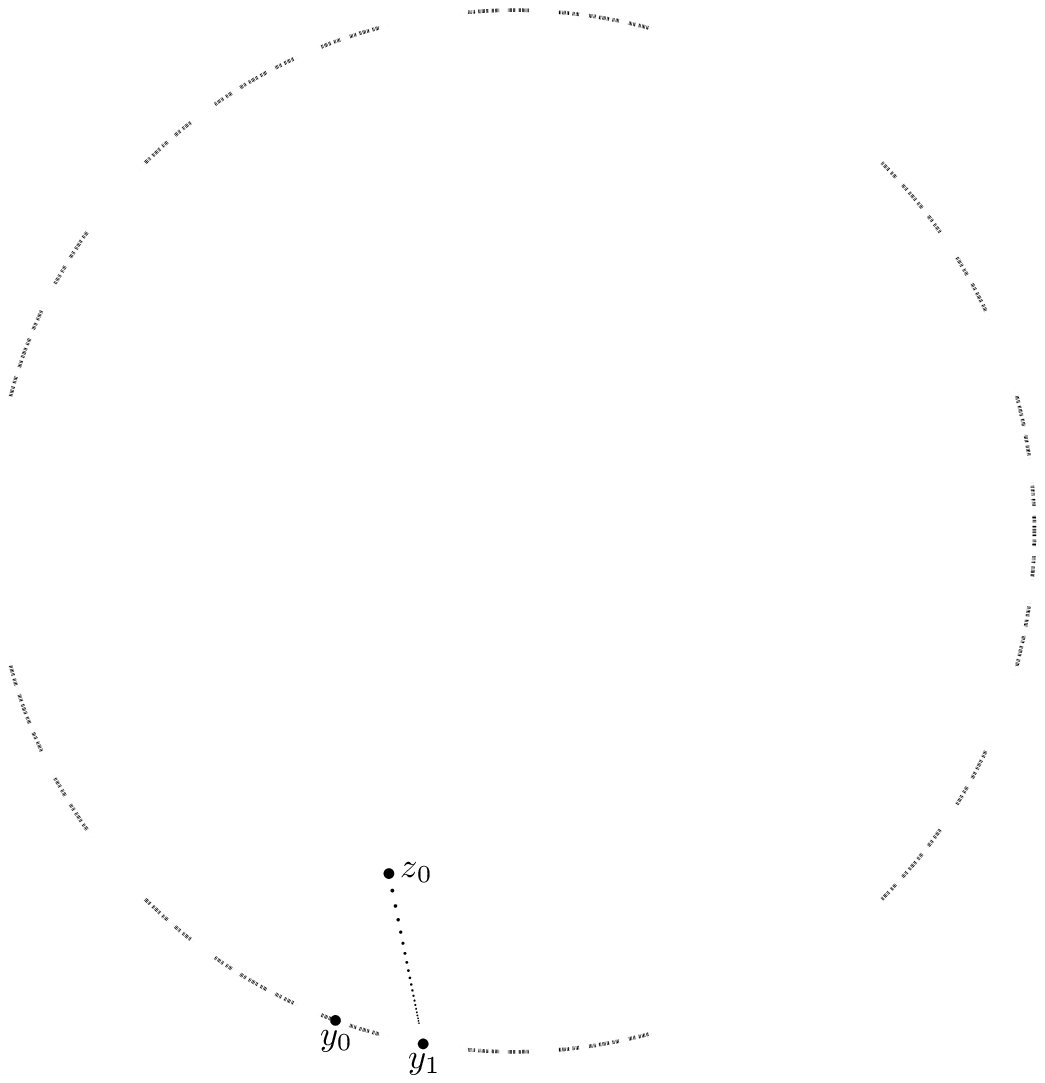}}
\caption{The terminal states are the unit sphere, and the support of $F$ is $S_\delta$, the $\delta$-neighborhood of a porous set on the unit sphere.  From a starting point $z_0$, player~\II finds a point $y_1$ on the sphere that is near the closest point $y_0$ on the sphere but far from $S_\delta$.  Player~\II then tugs towards $y_1$ along a sequence of points, and then when it gets much closer to the sphere than to $S_\delta$, it tugs straight to the sphere. \label{cantor}}
\end{figure}
It can be checked that after $k=O(-\log\alpha)$
steps the sequence hits $Y\setminus S$ and that the assumptions
of the lemma hold with this $k$, with $\gamma = \alpha/10$
and with some $\eps_0>0$, independent of $z_0$.
The theorem now easily follows from the lemma.
\end{proof}

We now proceed to study the $\infty$-harmonic measure of spherical caps.

\begin{proof}[Proof of Theorem~\ref{t.hmeasure}]
It turns out to be more convenient to work
with $-A_\delta=\{-y:y\in A_\delta\}$ in place of $A_\delta$.
By an obvious comparison argument,
it is sufficient to estimate $u(0)$ where
$u=u_\delta$ is the AM function in $X\setminus Y$
with boundary values
$$
F(y):=\begin{cases} 1\,,& y\in -A_\delta\,,\\
  0\,,& d(y,-A_\delta)>\delta\,,\\
  1-d(y,-A_\delta)\,\delta^{-1}\,, & 0<d(y,-A_\delta)\le \delta\,.
\end{cases}
$$
The function $u$ is invariant under rotations of $X$ preserving $(1,0,\dots,0)\in\R^n$.
Therefore, $u$ is also AM in $\R^2\cap X$.
Thus, we henceforth restrict to the case $n=2$, with no loss of generality.

Aronsson~\cite{MR850366} constructed a family of viscosity solutions $G_m$ to
$\Delta_\infty u=0$ in $\R^2\smallsetminus\{0\}$ that are separable in
polar coordinates: $G_m(r,\theta)=r^{m^2/(2m-1)} h_m(\theta)$ (for
$m\in\Z$).  We are interested in the $m=-1$ solution,
which may be written as
\begin{equation}\label{e.Aro}
G=\left[
\frac{\cos \theta\,(1-|\tan(\theta/2)|^{4/3})^2}{1 + |\tan(\theta/2)|^{4/3} +
  |\tan(\theta/2)|^{8/3}}\right]^{1/3} r^{-1/3} .
\end{equation}
Observe that when
$-\pi/2\leq\theta\leq\pi/2$ this solution is non-negative, and
$G\ge c\, r^{-1/3}$ when $\theta\in[-\pi/4,\pi/4]$, say,
where $c>0$ is some fixed constant.
Then $F\le O(1)\,\tilde G_\delta$, where
$$
\tilde G_\delta(x,y)= \delta^{1/3}\, G(x+1+2\,\delta,y)\,.
$$
Since $u(0)$ is monotone in $F$, and $\tilde G_\delta$ satisfies $\Delta_\infty\tilde G_\delta=0$,
it follows that $u(0)\le O(1)\,\tilde G_\delta(0)=O(\delta^{1/3})$.

We now show that the bound $u(0)=O(\delta^{1/3})$ is tight.
It is easy to see, using comparison with a cone centered
at $(-1,0)$, say, that $u(-1+\delta/10,0)>c>0$, where $c$
is a constant that does not depend on $\delta$.
Comparison with a cone centered at any point $z$
in the unit disk shows that
$u(z')\ge u(z)/2$ if $|z-z'|\le (1-|z|)/2$.
Using such estimates, it is easy to see that
there is a constant $c'>0$ such that
$u\ge c'$ on the disk $B$ of radius $\delta$
centered at $q:=(-1+2\delta,0)$, provided that $\delta<1/4$, say.
Now consider the function
$$
G^*_\delta(x,y)=c'\,\delta^{1/3}\,G\bigl((x,y)-q\bigr)\,.
$$
By the choice of the constant $c'$, we have $G_\delta^*\le u$ on $\partial B$,
since $G\le r^{-1/3}$ in $\R^2\setminus\{0\}$.
If $(r',\theta')$ denote the polar coordinates centered at $q$, the center of $B$,
and $(r,\theta)$ denote the standard polar coordinates centered at $0$,
then $2\,|\theta'|\ge \pi- |\theta-\pi|$ when $r=1$ and we choose $\theta'\in[-\pi,\pi)$.
Also, $r'$ is bounded from below by a constant times $|\theta-\pi|$
when $\theta\in [0,2\pi]$ and $r=1$.
Since $|G|\le O(1)\,\bigl((\pi/2)-|\theta|\bigr)\,r^{-1/3}$
when $\theta\in[-\pi,\pi)$, it follows that on the unit circle
$$
G_\delta^*\le O(1)\,\bigl((\pi/2)-|\theta'|\bigr) (r')^{-1/3}
\le
O(1)\,\delta^{1/3}\,|\theta-\pi|^{2/3}\le O(\delta^{1/3})\,.
$$
Therefore $G_\delta^*-O(\delta^{1/3})\le u$ on the boundary of the unit circle,
as well as on $\partial B$. Consequently, $u\ge G_\delta^*-O(\delta^{1/3})$
in the complement of $B$ in the unit disk.
This implies that there is some $r\in(0,1)$
such that $u_\delta(-r,0)\ge \delta^{-1/3}$ for all sufficiently small $\delta$.
The required estimate $u_\delta(0)\ge \Theta(\delta^{-1/3})$ now follows
by several applications of comparison with cones (the number of which
depends on $r$), similar to the above argument estimating
a lower bound for $u$ on $B$.
\end{proof}

One could also try to use Aronsson's other solutions $G_m$ to
bound $\infty$-harmonic measure of certain sets in other domains.
Also, Aronsson has some explicit solutions for $\Delta_p u=0$,
which may serve a similar purpose.

\section{Counterexamples} \label{counterexamplesection}

\subsection{Tug-of-war games with positive payoffs and no value}
\label{sub:pospay}

Here we give the promised counterexample showing that
the hypothesis $\inf f >0$
in Theorem~\ref{towvalue}
 cannot be relaxed to $f> 0$.
As we later point out,
a similar construction works in the continuum setting
of length spaces.

\begin{figure}[hb]
\centerline{\includegraphics[height=2in,bb=97 97 188 278]{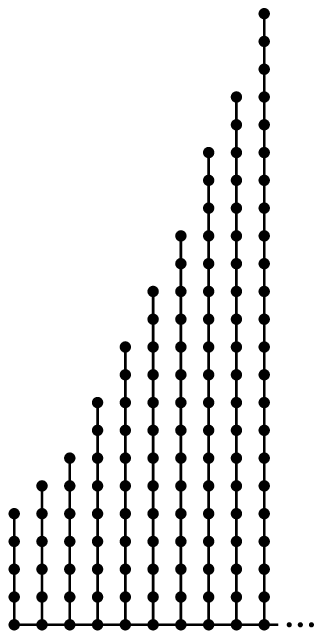}}
\caption {\label{f.comb} A comb.}
\end{figure}

\indent
The comb game is tug-of-war with running payoffs on an infinite
graph shaped like a comb, as shown in Figure~\ref{f.comb}.
The comb is defined by an infinite sequence of
positive integers $\ell_0,\ell_1,\ell_2\ldots$, and has states (vertices)
$\{(x,y)\in\Z^2:0\leq x\text{ and } 0\leq y \leq \ell_x\}$.  The edges of
the comb are of the form $(x,y)\sim(x,y+1)$ and
$(x,0)\sim(x+1,0)$, giving the graph the shape of a comb, where the
$x\thh$ tooth of the comb has length $\ell_x$.
The running payoff $f(x,y)$ is $1/\ell_x$ if $y=0$ and zero otherwise.
(Later we will consider a variation where $f>0$ everywhere.)
The terminal states are the
states of the form $(x,\ell_x)$, and the terminal payoff $F$
is zero on all terminal states.

\begin{lemma}\label{l.pullup}
  For any comb (choice of the sequence $\{\ell_x\}$), we have $u_\I(0,0)=2$.
\end{lemma}
\begin{proof}
First we argue that $u_\I(0,0) \geq 2$.
Denote by   $\psi(t):=\sum_{i=0}^{t-1} f(x_{i},y_{i})$ the accumulated running payoff.
Fix a large $B>0$ and equip player~\I with the following strategy:
 at all times $t < \tau$ for which $\psi(t) < B$,
player~\I pulls down if $y_t \not = 0$, left if $y_t = 0$ and $x_t
\not = 0$, and right if $(x_t, y_t) = (0,0)$;
if $\psi(x_t) \ge B$, then player~\I pulls toward the closest terminal state.
Define the termination time $\tau=\inf\{t: y_t=\ell_{x_t}\}$
and also the stopping time $\sigma_0= \inf \{t: x_t=0 \}\wedge \tau$.
Fix any strategy for player~\II and observe that  $\{x_{t\wedge \sigma_0} \}_{t \ge 0}$ is a nonnegative
supermartingale, which must converge a.s.; therefore, from any initial state,
$\sigma_0<\infty$ a.s.; therefore, the stopping time $\sigma = \inf \{t: \psi(t) \geq B \}\wedge \tau$
is almost surely finite, whence $\tau<\infty$ a.s.\ as well.

 Consider the process
  $$M_t = 2(1-y_t/\ell_{x_t}) + \psi(t) \,.$$
  Then $M_{t\wedge \sigma}$ is a submartingale.
Note that if
$y_\sigma = 0$ then $M_\sigma = \psi(\sigma) + 2$ and $\psi(\sigma) \geq B$, while
if $y_\sigma >0$ then $\sigma=\tau$ and $M_\sigma = \psi(\sigma)$.  In any case, $M_\sigma
\leq \psi(\sigma)(1 + 2/B)$. Thus by optional stopping, $$2 = M_0 \leq
\mathbb E M_\sigma \leq (1 + \frac{2}{B})\mathbb E \psi(\sigma) \leq (1
+ \frac{2}{B}) u_\I(0,0).$$

  Next, to show that $u_\I(0,0) \leq 2$,
  suppose that player~\II adopts the strategy of always pulling up, and
  player~\I knows this and seeks to maximize her payoff.
  Because player~\II always pulls up, $M_t$ is a positive
  supermartingale, so $M_t$ a.s.\ converges to $M_\infty$, and by
  optional stopping, $$2(1-y_0/\ell_{x_0})=M_0\geq
  \E[M_\infty] \geq \E\big[\sum_{t'=0}^\infty
  f(x_{t'},y_{t'})\big].$$
  This shows that by always pulling up, player~\II can force player~\I's
  expected payoff from $(x_0,y_0)=(0,0)$ to be no more than $2(1-y_0/\ell_{x_0})=2$.
\end{proof}
\begin{remark}
 Note
  that the strategy for player~\II of always pulling up does not necessarily force the
  game to end with probability one;  this strategy always gives an upper
  bound on $u_\I(0,0)$, but it only yields an upper bound on $u_\II(0,0)$ if $\sum_{x \ge 0} \ell_x^{-1}=\infty$, when the Borel-Cantelli Lemma ensures
  termination.
\end{remark}

Suppose that the teeth of the comb are long, i.e.\
$\sum_{x=0}^\infty \ell_x^{-1} <\infty$, and that player~\I plays the
strategy of always pulling down if $y_t>0$ and always pulling right if $y_t=0$.
If player~\II plays the strategy of always pulling up,
then we may calculate the probability that the game terminates when
started in state $(x,0)$.  Either the state goes to the terminal
state of the tooth or goes to the base of the next tooth, and the
latter probability is $\ell_x/(\ell_x+1)$.  When the teeth of the
comb are long, the game lasts forever with probability $\prod_{x \ge 0}\ell_x/(\ell_x+1)>0$.  If
player~\II needs to ensure that the game terminates, she will need
to be prepared to sometimes pull left instead of up,
and this necessity will be costly for player~\II.

\begin{lemma}\label{l.pullright}
For every $s>0$ there is a suitable comb (choice of the sequence
$\{\ell_x\}$) such that $u_\II(0,0) \geq s$.
\end{lemma}

\begin{proof}
Let player~\I play the strategy of pulling down when $y>0$ and
pulling to the right when $y=0$.
We will estimate from below the expected payoff when player~\II
plays any strategy
which guarantees that the game terminates in finite
time a.s.\ against the above strategy for player~\I.
Let $a_x$ denote the probability that the game
terminates at the terminal state $(x,\ell_x)$.
Then $ \sum_{x=0}^\infty a_x=1$.
For every $x\in\N$ let $L_x,U_x$ and $R_x$
denote the expected number of game
transitions from $(x,0)$ to the left, upwards and to the right,
respectively (that is, to $(x-1,0)$, to $(x,1)$ and to $(x+1,0)$,
respectively).

First, observe that
$$
U_x\ge a_x\,\ell_x\,,
$$
because every time that the game state arrives
at $(x,1)$, with conditional probability at most $1/\ell_x$ the
game terminates at $(x,\ell_x)$ before returning
to $(x,0)$. Next, observe that
$$
L_{x+1}\ge R_x-1\,,
$$
because the number of transitions from
$(x,0)$ to $(x+1,0)$ can exceed the number of
transitions from $(x+1,0)$ to $(x,0)$ by at most
one. (More precisely, we have $L_{x+1}= R_x-\sum_{j>x} a_j$,
but this will not be needed.)
Finally, note that
$$
R_x\ge L_x+U_x\,,
$$
because player~\I always pulls right at $(x,0)$.
(We set $L_0:=0$.)

An easy induction shows that the above
relations imply
$$
R_x\ge \sum_{j=0}^x a_j\,\ell_j - x
\,.
$$
The expected payoff $\mathcal {E}$ satisfies
$$
\mathcal E=\sum_{x=0}^\infty \frac{R_x+L_x+U_x}{\ell_x}
\ge
\sum_{x=0}^\infty\frac{ R_x}{\ell_x}\,.
$$
We plug in the above lower bound on $R_x$
to obtain
$$
\mathcal E \ge
\sum_{j=0}^\infty a_j\,\Bigl(\sum_{x=j}^\infty \frac{\ell_j}{\ell_x}
\Bigr)
-\sum_{x=0}^\infty \frac x{\ell_x}\,.
$$
Since $\sum_{x\ge 0}a_x=1$, the lemma follows if we choose $\ell_x=(c+x)^3$
with $c>0$ sufficiently large.
\end{proof}

Note that when player~\II always pulls up, for every vertex $v$ in
the comb there is a finite upper bound $b(v)$ on the expected number
of visits to $v$, which holds regardless of the strategy used by
player~\I. Since in the proof of Lemma~\ref{l.pullup} player~\II
always pulls up, it follows that the value for player~\I is at most
$3$ even when $f$ is replaced by $f(\cdot)+
q(\cdot)/\bigl(b(\cdot)+1\bigr)$, where $q>0$ and $\sum_v q(v)=1$.
This modification will certainly not decrease the value for
player~\II. Therefore, in Theorem~\ref{towvalue}, the assumption
$\inf f>0$ cannot be replaced by the assumption $f>0$, even if $F=0$
throughout $Y$.

Note that we may convert the discrete comb graph to a continuous length
space by adding line segments corresponding to
edges in the comb.
It is then easy to define the corresponding continuous $f$
and conclude that
in Theorem~\ref{continuumvalueexists} the assumption $\inf |f|>0$ cannot be relaxed to
$f>0$.
The details are left to the reader.

\begin{remark}
The dependence of $u_\II(0,0)$ on the growth of $\{\ell_x\}$ is somewhat
surprising.  If $\{\ell_x\}$ grows slowly enough so that $\sum_{x \ge 0}
\ell_x^{-1} = \infty$, then $u_\II(0,0) \leq 2$, because (as remarked above) the
Borel-Cantelli lemma implies that the strategy of always pulling up
is guaranteed to terminate, and thus the proof of Lemma
\ref{l.pullup} applies.  We have seen that for some sequences $\{\ell_x\}$ of
polynomial growth, $u_\II(0,0)$ can be arbitrarily large.  However,
if $\{\ell_x\}$ grows rapidly enough so that $\ell_{x+1}/\ell_x \to
\infty$, then $u_\II(0,0) \leq 2$.  This is immediate from the
following more general statement: for every $k
> 0$, we have $u_\II(0,0) \leq 2+2 \ell_k \sum_{j > k} \ell_j^{-1}.$  To
prove this, suppose player~\II adopts the strategy of always pulling
left when $x_t > k$ and $y_t = 0$ and up otherwise.  Let $\psi(t,k)$
denote the payoff accumulated for player~\I at points in $[0,k]
\times \{0\}$ up to time $t$.  Then $\mathbb E \psi(t,k) \leq 2$,
because
$$\widetilde M_t = 2(1 - \frac{y_t}{\ell_{x_t}}) + \psi(t,k)$$
is a positive supermartingale with $\widetilde M_0=2$.  Then we claim that $$\mathbb E[
\psi(t) - \psi(t,k)] \leq 2 \ell_k \sum_{j>k} \ell_j^{-1}.$$

Fix $j>k$. Starting at $(k,0)$, the expected number of visits
to $(j,0)$ before returning to $(k,0)$ is at most 1 (by comparison to simple random walk).
The expected number of visits to $(k,0)$ is at most $2\ell_k$, since
each time $(x_t,y_t) = (k,0)$, there is a chance of at least $1/(2\ell_k)$ of
terminating at $(k,\ell_k)$ without returning to $(k,0)$. Thus the expected accumulated payoff
at $(j,0)$ is at most $ 2 \ell_k \ell_j^{-1}$; summing over $j>k$ proves the claim.
\end{remark}

\subsection{Positive Lipschitz function with multiple AM extensions}
\label{multiple-AM}

Here we show that uniqueness in Theorem~\ref{uniqueAM} may fail if $F$
is unbounded, that is, we give an $(X,Y,F)$ (here $f=0$) for which $F$
is Lipschitz and positive and the continuum value is not the only AM
extension of $F$.

Let $T$ be the rooted ternary tree, where each node has three direct
descendants and every node but the root has one parent.  Let $X$ be
the corresponding length space, where we glue in a line segment of
length $1$ for every edge in $T$. For every node $v$ in $T$, we
label the three edges leading to descendants of $v$ by $1$, $-1$ and
$*$.  (One may interpret $T$ as the set of finite stacks of cards,
where each card has one of the three labels $1$, $-1$, and $*$.)

We define a function $w$ on the nodes of $T$ by induction on the
distance from the root.
For any vertex $v$ let $k(v)$ be the number of edges
on the simple path from the root to $v$ whose label
is not $*$.
Set $w(\text{root}):=0$.
Next, if $v$ is the parent of $v'$ and the edge from
$v$ to $v'$ is labeled $\pm 1$, then $w(v')=w(v)\pm 1$,
respectively.
Finally, if the edge from $v$ to $v'$ is labeled $*$, let
$w(v')=w(v)+1-2^{-k(v)}$, say.
This defines $w$ on the vertices of $T$. We define it
on $X$ by linear interpolation along the edges.

Let $V_0$ denote the set of nodes
consisting of the root and all vertices
$v$ such that  the edge from $v$ to its parent
is not labeled $*$.
For each $v\in V_0$ let $b(v)$ be some large integer,
whose value will be later specified, and let
$q(v)$ be the vertex at distance $b(v)$ away
from $v$ along the (unique) infinite simple path starting
at $v$ which contains only edges labeled $*$.
Let $Y_1=\{q(v):v\in V_0\}$ and
let $Y_0$ be the set of all nodes $v$
such that $w(v)\le 3/2$. Set $Y:=Y_0\cup Y_1$.

We first claim that $w$ is AM on $X\setminus Y$.
Indeed, let $U$ be an open subset of $X\setminus Y$.
If $U$ does not contain any tree node, then $\Lip_Uw=\Lip_{\partial U}w$,
because $w$ interpolates linearly inside the edges.
Suppose now that $v\in U$ is a node.
Let $\beta$ be the infinite path starting at $v$ going always
away from the root which uses
only edges labeled $-1$.
For each positive integer $m$ let $\gamma_m$ be the infinite path
starting at $v$ going always away from the root whose
first $m$ edges are labeled $1$ and the rest are labeled $*$.
Observe that $\gamma_m$ meets $Y_1$ and
$\beta$ meets $Y_0$. Consequently, $\gamma_m\cap\partial U\ne\emptyset$
and $\beta\cap\partial U\ne\emptyset$.
If $x_1\in \gamma_m\cap\partial U$ and $x_0\in\beta\cap\partial U$,
then $w(x_1)-w(x_0)\ge (1-2^{-m})\,d(x_1,x_0)$ by the construction of $w$.
Thus $\Lip_{\partial U} w\ge 1$. Since $\Lip_X w=1$, this proves that
$w$ is AM on $X\setminus Y$.

Next, we consider the discrete tug-of-war game on the graph $T$
where $f=0$ and $F$ is the restriction of $w$ to $Y$.  Let $x_0$, the
starting position of the game, be the third vertex on the infinite
simple path from the root whose edges are all labeled $1$.  We claim
that the value $u_\I$ for player~\I satisfies $u_\I(x_0)< w(x_0)$.
Before proving this we explain the idea: at each step, player~\I has
one or two moves that increase the value of $w$ by $1$ (either
moving away from the root along an edge with label $1$ --- i.e.,
``adding a $1$ card to the deck'' --- or moving towards the root
along an edge labeled $-1$ --- i.e., ``removing a $-1$ card from the
top of the deck''), and similarly, player~\II has one or two moves
that decrease the value of $w$ by $1$. Player~\I can thus make
$w(x_j)$ a submartingale by always making moves of this type.

This does not force the game to terminate, however; in order to end
the game favorably by reaching a point in $Y_1$, player~\I must add
a sequence of cards labeled $*$ to the top of the deck.  If
player~\II adopts the strategy of always choosing the edge labeled
$-1$ (``adding a $-1$ card to the deck''), then (provided that
$b(\cdot)$ increases rapidly enough) the expected number of
suboptimal moves that player~\I must make (by moving along an edge
with label $*$) to reach any particular element of $Y_1$ is large
enough to significantly decrease the total expected payoff for
player~\I.

We proceed to prove that $u_\I(x_0)<w(x_0)$. Let $v_0$ be a vertex
in $V_0\setminus Y$, which is in the same connected component of
$T\setminus Y$ as $x_0$. Let $(v_0,v_1,\dots,v_{b(v_0)})$ be the
simple path from $v_0$ to $q(v_0)$. We now abbreviate $n=b(v_0)$,
$h=w(v_0)$ and $a= 1-2^{-k(v_0)}$.

Let player~\II use the naive strategy of
always trying to move from the current state along
the edge labeled $-1$ going away from the root.
For $i\in\{0,1,\dots,n\}$ let
$\tau_i$ be the first time $t$ such that $x_t=v_i$,
and if no such $t$ exists let $\tau_i=\infty$.
(As before $x_t$ is the game position at time $t$.)
Let $M^{(i)}_t = w(x_t)$ for $t<\tau_i$
and $M^{(i)}_t= w(v_{i-1})+1$ for $t\ge\tau_i$.
Then $M^{(i)}_t$ is clearly a supermartingale, regardless
of the strategy used by player~\I.
Since the game terminates at time $t$ if $w(x_t)\le 1$,
this implies that for $i\in \N_+$
$$
\Pr\bigl[ \tau_{i}<\infty \bigm| \tau_{i-1}<\infty\bigr]
\le \frac{w(v_{i-1})}{w(v_{i-1})+1} =
1-
\frac {1}{h+(i-1)\,a +1}
\le \exp\Bigl(
\frac {-1}{h+(i-1)\,a +1}\Bigr)
\,.
$$
Consequently,
\begin{multline*}
\Pr\bigl[ \tau_{n}<\infty \bigr]
\le
\prod_{i=1}^n
\Pr\bigl[ \tau_{i}<\infty \bigm| \tau_{i-1}<\infty\bigr]
\le
\exp\Bigl( -\sum_{i=1}^{n}
\frac {1}{h+(i-1)\,a +1}\Bigr)
\\
\le
\exp \Bigl(-\int_{0}^n \frac {ds}{h+1+s\,a}\Bigr)
=
\Bigl( \frac{h+1}{a\,n+h+1}\Bigr)^{1/a}
\,.
\end{multline*}
 Therefore
$$
\Pr\bigl[{\tau_n<\infty} \bigr]F(v_n)
\le (h+n\,a)\,
\Bigl( \frac{h+1}{a\,n+h+1}\Bigr)^{1/a}.
$$
Since $a<1$, we can make this smaller than any required positive
number by choosing $n=b(v_0)$ sufficiently large.
We may therefore choose the function $b:V_0\to\N_+$ so
that
$$
\sum_{v\in Y_1}
\Pr\bigl[ \text{game ends at $v$}\bigr] \,F(v)
<1/2\,.
$$
Since $F\le 3/2$ on $Y_0$, this implies $u_\I(x_0)<2<3=w(x_0)$.

Let $u_*$ be the linear interpolation of $u_\I$ to the edges.
By Lemma~\ref{valuesareinfharmonic} $u_\I$ is discrete $\infty$-harmonic.
It follows that in the \II-favored $\eps$-tug-of-war game, for every $\delta>0$
and $\eps\in (0,1/2)$
player~\II can play to make $u_*(x_t)+2^{-t}\delta$
a supermartingale. Consequently,
Lemma~\ref{vugap} implies that the value of the continuum
game is bounded by $u_\I$ on the vertices.
Since the continuum value is AM on $X\setminus Y$,
by Theorem~\ref{uniqueAM}, this proves our claim that
the assumption $\sup F<\infty$ is necessary for uniqueness
to follow in the setting of that theorem.

\subsection{Smooth $f$ on a disk with no game value}
\label{sub:noval}

The following example is a continuum analog of the triangle example
at the beginning of Section \ref{s.run}.

Let $X$ be a closed disk in $\R^2$, and let $Y=\partial X$
be its boundary.
We now show that for some smooth function $f:X\to\R$,
 the value $u^\eps_\I$ for player~\I and the value
$u^\eps_\II$ for player~\II differ in $\eps$-tug-of-war with
$F=0$ on $Y$. Not only are the values different,
but the difference does not shrink to zero as $\eps\searrow 0$.
There is a lot of freedom in the choice of the function $f$,
but an essential property, at least for the proof, is that $f$
is anti-symmetric about a line of symmetry of $X$.

It will be convenient to identify $\R^2$ with $\C$, and
use complex numbers to denote points in $\R^2$.
Let $f$ be a $C^\infty$ function such that:
\begin{enumerate}
\item $f\ge 1$ inside the disk $|z-1|<1/2$,
\item $0\le f\le 2$ in the right half plane $\Re z\ge 0$,
\item $f=0$ in $\{z\in\C: \Re z\ge 0,\, |z-1|>11/20\}$, and
\item $f$ is anti-symmetric about the imaginary line.
\end{enumerate}
Let $R>2$ be large, let $X=\{z:|z|\le R\}$, $Y=\{z:|z|=R\}$
and $F=0$ on $Y$.
Consider $\eps>0$ small.
We want to show that for an appropriate choice of $R$, then
\begin{equation}\label{e.noval}
\liminf_{\eps\searrow 0} \|u_\I^\eps-u_\II^\eps\|_\infty >0\,.
\end{equation}
Indeed, suppose that this is not the case.
Let $\delta>0$ be small, and suppose that
$\|u_\I^\eps-u_\II^\eps\|_\infty <\delta$.
Symmetry gives $u_\I^\eps (x+i\,y)=-u_\II^\eps(-x+i\,y)$.
The assumption
$\|u_\I^\eps-u_\II^\eps\|_\infty <\delta$ therefore tells us that
$|u_\I^\eps(x+i\,y) + u_\I^\eps(-x+i\,y)|< \delta$.
In particular,
$|u_\I^\eps|< 2\,\delta$ on the imaginary axis.
We now abbreviate $w=u_\I^\eps$.
Clearly $w\ge -\delta$ on $\{z\in X: \Re z\ge 0\}$
and $w\le \delta$ on $\{z\in X:\Re z\le 0\}$.
By considering strategies which pull towards the imaginary
axis and then using whatever strategy gives
a value in $[-2\,\delta,2\,\delta]$ it is easy to
see that  $|w| = O(1)$
(as usual, $f=O(g)$ means that there is some universal
constant $C$ such that $f\le C\,g$)
and that $|w|\le O(\delta+\eps)$
in $\{z\in X: |\Re z|\le 2\,\eps\}$, say.
(See, e.g., the proof of Lemma~\ref{uniformlylipschitz}.)
It is then easy to see that
there is a constant $c_1>0$, which does not depend on $R$ and
$\eps$ (as long as $\eps$ is sufficiently small and $R>2$, say)
such that $w \ge c_1$ on the disk $|z-1|<3/4$
and $w\le -c_1$ on the disk $|z+1|<3/4$.
Since $w$ is nearly anti-symmetric, it is enough to prove
the first claim.
Indeed, a strategy for player~\I which demonstrates this
is one in which she pulls towards $1$ until she accumulates
a payoff of $1$. If successful, she then aims towards the
boundary $Y$, but whether successful or not,
whenever the game position comes within distance $\eps$ of
the imaginary axis, she changes strategy and adopts an arbitrary
strategy that yields a payoff $O(\delta+\eps)$.
We assume that $\delta$ and $\eps$ are sufficiently small
so that the $O(\delta+\eps)$ term is much smaller than $c_1$.

It now easily follows from Theorem~\ref{t.hmeasure}
that $|w(z)|\le O(1)\,|z|^{-1/3} + O(\delta+\eps)$.
Consequently, there is some constant $r_1>2$ such that
$|w|< c_1/10$ on
$\{z\in X:|z|\ge r_1\}$, provided that $\eps+\delta$
is sufficiently small.
Set $L=c_1/(8\,r_1)$, and let $Z_L=Z_L^\eps$ be the
set of points $z\in X$ such that
$$
\sup \{w(z'):z'\in B_z(\eps)\}-
\inf \{w(z'):z'\in B_z(\eps)\} > 2\,L\,\eps\,.
$$
Since $|w|\le c_1/10$ on $|z|=r_1$
and $|w|\ge c_1$ on $S:=\{|z-1|<3/4\}\cup\{|z+1|<3/4\}$,
it follows that every path in the graph $(X,E_\eps)$ from $S$ to
$Y$ must intersect $Z_L$.
(We are assuming $R>r_1$.)
Let $\tilde Z_L$ denote the union of $Z_L$
and all connected components of $(X\setminus Z_L,E_\eps)$
that are disjoint from $Y$.

We now claim that there is some $r_2>r_1$ (which does not depend on
$\eps$) such that $Z_L\subset \{|z|<r_2\}$ and
therefore $\tilde Z_L\subset \{|z|<r_2\}$.
Indeed, if $z_1,z_2\in X$ satisfy $|z_1-z_2|<\eps$
and $|z_1|>3\,r_1$, say, then the strategies
for either player of pulling towards $z_2$
and only giving up when the current position is
within distance $\eps$ of $|z|=r_1$ show that
$$
|w(z_1)-w(z_2)|\le
\frac{ \eps}{|z_1|-r_1-2\,\eps}
\,\sup\{|w(z)|:|z|\ge r_1\}
\le
\frac{ \eps}{|z_1|-r_1-2\,\eps}\, c_1/10
\,.
$$
This proves the existence of such an $r_2$.

We now let  player~\II play a strategy very similar to the
backtracking strategy she used in the proof of Theorem~\ref{uniqueinfharmongraph}, which
we presently describe in detail.
Let $x_0,x_1,\dots$ denote the sequence of positions of the game.
Assume that $x_0\in \tilde Z_L$. For each $t\in\N$, let
$j_t:=\max\{j\in\N: j\le t,\, x_j\in \tilde Z_L\}$.
Note that if $j_t<t$, then $x_{j_t}\in Z_L$.
By Lemma~\ref{valuesareinfharmonic} inside
$\{z\in X: |z-1|\ge 3/4,\, |z+1|\ge 3/4\}$
the function $w$ satisfies $\Delta_\infty w=0$.
If $x_t\in \tilde Z_L$ and player~\II gets the turn,
she moves to some $z\in B_{x_t}(\eps)$ that
satisfies $w(z)< \inf\{w(z'):z\in B_{x_t}(\eps)\}+\delta\,2^{-t-1}$.
While if $x_t\notin \tilde Z_L$ and player~\II gets the turn,
she moves to any neighbor $z$ of $x_t$
that is closer (in the graph metric) to $x_{j_t}$ than $x_t$ in the subgraph
$G_t$
of $(X,E_\eps)$ spanned by the vertices $x_{j_t},x_{j_t+1},\dots,x_t$.

Now consider the game evolving under any strategy for player~\I and
the above strategy for player~\II. Let $d_t$ be the graph-distance
from $x_t$ to $x_{j_t}$ in $G_t$. Set
$$
m_t := w(x_{j_t})+ d_t\,\eps\,L+\delta\,2^{-t}+\sum_{k=0}^t f(x_k)\,.
$$
As in the proof of Theorem~\ref{uniqueinfharmongraph},
it is easy to verify that
$w(x_t)\le m_t$ and that $m_t$ is a supermartingale.

Let $\tau$ be the time in which the game stops; that is, $\tau:=\inf\{t\in\N :x_t\in Y\}$.
Assume that $\tau<\infty$ a.s.
Suppose that we could apply the optional stopping time theorem.
Then we would have
$$
w(x_0) \ge \E\bigl[m_\tau\bigr]=\E[ w(x_{j_\tau})+\delta\,2^{-\tau}] +\eps\,L\, \E[d_\tau] + \E\bigl[\text{payoff}\bigr].
$$
Thus,
$$
\E\bigl[\text{payoff}\bigr]\le O(1)-\eps\,L\,\E[d_\tau]\,.
$$
But $d_\tau\ge (R-r_2)/\eps$, because the Euclidean distance from $Y$ to $\tilde Z_L$
is at least $R-r_2$.
Since $w(x_0)$ is at most the supremum of expected payoffs under the current strategy for
player~\II and an arbitrary strategy guaranteeing $\tau<\infty$ for player~\I,
we obtain $ w(x_0)\le  O(1)- R+r_2$. This contradicts our previous conclusion that $\|w\|_\infty= O(1)$,
because we may choose $r$ large.

How do we justify the application of the optional stopping time theorem?
We slightly modify the strategy for player~\II.
So far, our analysis utilized one advantage to player~\II in the calculation of $w=u_\I^\eps$; that is,
that it is player~\I's responsibility to make sure that $\tau<\infty$.
Now we have to use the other advantage to player~\II, which is that player~\II gets to choose
her strategy after knowing what player~\I's strategy is.
Let $t_0$ be the first time such that with the strategy which player~\I is using and the above
strategy for player~\II, the game terminates by time $t_0$ with probability at least $1/2$.
The new strategy for player~\II is to play the above strategy until time $t_0$, and if the
game lasts longer to use an arbitrary strategy which would guarantee
a conditioned expected future payoff of $u_{\II}^\eps(x_{t_0})=O(1)$.
We may certainly apply the optional stopping time theorem at time $\tau\wedge t_0$.
This does give a contradiction as above, because $\E[d_{\tau\wedge t_0}]\ge (R-r_2)/(2\,\eps)$,
and the contradiction proves our claim~\eqref{e.noval}.

\subsection{Lipschitz $g$ on unit disk with multiple solutions to $\Delta_\infty u = g$}

Here we construct the counterexample showing that if we omit the
assumption $\inf f>0$, in Corollary~\ref{c.linf}, then it may fail.
In the example, $X$ is the unit disk in $\R^2$, $Y$ is its boundary,
$F=0$ and $f$ is Lipschitz and take values of both signs. The
example is motivated by and similar to the example of
Section~\ref{sub:noval}. However, since the assumptions of
Corollary~\ref{c.linf} do not hold (obviously), we need to construct
the solutions $u$ to $\Delta_\infty u=g$ not by using tug-of-war,
but by other means. In fact, we will use a smoothing of Aronsson's
function~\eqref{e.Aro}.

The following analytic-geometric lemma will replace
the use of the backtracking strategy for player~\II in Section~\ref{sub:noval}.
The following notation will be needed. For $x\in U$ and $u:R\to\R$ let
$$
\Lip_xu:= \inf\{\Lip_W u:x\in W\subset U,\,W\text{ open}\}\,.
$$

\begin{lemma}\label{l.nocontin}
Let $X$ be a length space, and let $F_1,F_2$ be bounded Lipschitz real valued functions
on a nonempty subset $Y\subset X$. Let $u_1$ and $u_2$ be the corresponding AM
extensions to $X\setminus Y$.
Fix $L>0$, and let $Z_L:=\{x\in X:\Lip_x u_1> L\}$.
Suppose that $F_1=F_2$ on $Z_L\cap Y$, that $\Lip_{Y\setminus Z_L} F_2 \le L$ and
\begin{equation}\label{e.nomatter}
\sup\Bigl\{\frac{|F_2(y)-u_1(x)|}{d(y,x)}:{y\in Y\setminus Z_L},\,{x\in Z_L}\Bigr\} \le L\,.
\end{equation}
Then $u_1=u_2$ in $Z_L$.
\end{lemma}

We need a few simple observations before we begin the proof.
Suppose that $U$ is an open connected subset of a length space $X$,
that $X\setminus U$ is nonempty,
and that $F:X\setminus U\to\R$ is Lipschitz and bounded.
Let $\tilde\partial U$ denote the set of all points
$x\in\partial U$ such that there is a finite length path $\gamma\subset X$
where $\gamma\cap U\ne\emptyset$ and $\gamma\setminus U=\{x\}$.
Set $U_*:=U\cup\tilde\partial U$
and for two points
$x,y\in U_*$ let $d_*(x,y)=d_*^U(x,y)$ be the infimum
length of paths $\gamma\subset U_*$ joining $x$ and $y$
such that $\gamma\cap \tilde\partial U$ is finite.
Note that $(U_*,d_*)$ is a length space.
Also, since $d_*\ge d$ on $U_*\times U_*$, the restriction
of $F$ to $\tilde\partial U$ is Lipschitz.
Consequently, this restriction of $F$
has an AM extension $u:U_*\to\R$ with respect to $d_*$.
We claim that this is also an AM extension of $F$ with respect to $d$.

Observe that for any Lipschitz function $w:V\to\R$ where
$V\subset X$ is open, we have
$$
\Lip _{V} w \le \Bigl(\sup_{x\in V} \Lip_xw\Bigr)\vee\bigl( \Lip_{\partial V} w\bigr)\,.
$$
(This can be verified by considering for every pair of distinct points $x,y\in V$
the restriction of $w$ to a nearly shortest path joining $x$ and $y$ in $X$.)
Hence, to show that $u$ is AM with respect to $d$ it is enough to prove
that
\begin{equation}\label{e.inter}
\sup_{x\in V} \Lip_xw\le  \Lip_{\partial V} w
\end{equation}
holds for arbitrary open $V\subset U$.
Observe that for $x\in U$ we have
$\Lip_x w=\Lip^*_x w$, where $\Lip^*$ refers to $d_*$.
Since $w$ is AM with respect to $d_*$, we have
$\Lip^*_x w\le \Lip^*_{\partial V}w$ for all $x\in V$.
Finally, $\Lip^*_{\partial V}w\le \Lip_{\partial V}w$ because $d_*\ge d$ on
$U_*\times U_*$. This proves~\eqref{e.inter} and thereby shows that $w$ is
AM on $U$ with respect to $d$. In particular, we conclude that
\begin{equation}\label{e.upshot}
\Lip^*_{U_*} w =\Lip^*_{\tilde\partial U} w
\end{equation}
for the AM extension of $F$ to $U$.

\begin{proof}[Proof of Lemma~\ref{l.nocontin}]
First, it is clear that $Z_L$ is a closed set.
Let
$$
F_*(x)=\begin{cases} u_1(x)\,,& x\in Z_L\,,\\
  F_2(x)\,,& x\in Y\setminus Z_L\,,
\end{cases}
$$
and let $w:X\to\R$ be the AM extension of $F_*$ to $X\setminus (Y\cup Z_L)$.
We claim that $w$ is also an AM extension of $F_2:Y\to\R$ to $X\setminus Y$.
Since uniqueness of AM extensions holds in this setting,
this will imply $w=u_2$ and complete the proof, because
$w=u_1$ on $Z_L$.

Let $U$ be a connected component of $X\setminus Z_L$,
and let $d_*=d_*^U$ denote the corresponding metric
on $U\cup \tilde\partial U$.
Note that $\Lip_{\tilde\partial U}^*u_1\le L$,
since $\Lip_x u_1\le L$ for $x\in X\setminus Z_L$.
By~\eqref{e.nomatter} and the assumption
 $\Lip_{Y\setminus Z_L} F_2 \le L$ it follows that
that $\Lip_{\tilde \partial U\cup (Y\cap U)}^* F_*\le L$.
Thus, we get $\Lip_x w=\Lip_x^* w\le L$ for $x\in U$,
which implies $\Lip_x w\le L$ for $x\in X\setminus Z_L$.

Now let $V\subset X\setminus Y$ be open.
In order to prove that $w$ is AM, we need to establish~\eqref{e.inter}.
Without loss of generality, we assume that $V$ is connected.
If $V\cap Z_L=\emptyset$, then~\eqref{e.inter} certainly holds,
since $w$ is AM in $X\setminus (Z_L\cup Y)$.

Suppose now that $V\cap Z_L\ne\emptyset$.
Then $L_1:=\Lip_V u_1>L$, by the definition of
$Z_L$. Since $u_1$ is AM,
by~\eqref{e.upshot} with $u_1$ in place of $w$,
$V$ in place of $U$ and $d_*^V$ in place of $d_*^U$,
there is a sequence of pairs $(x_j,y_j)$ in $\tilde\partial V$
and a sequence of paths $\gamma_j\subset V\cup\tilde\partial V$
from $x_j$ to $y_j$, respectively,
such that $x_j\ne y_j$,
\begin{equation}\label{e.tight}
\lim_{j\to\infty}\frac{|u_1(x_j)-u_1(y_j)|}{\operatorname{length}(\gamma_j)}= L_1\,,
\end{equation}
and $\gamma_j\cap\tilde\partial V$ is finite.
For all sufficiently large $j$,
$|u_1(x_j)-u_1(y_j)|> L\,\operatorname{length}(\gamma_j)$.
Hence $\gamma_j\cap Z_L\ne\emptyset$.
Let $x_j'$ be the first point on $\gamma_j$ that is in $Z_L$ and
let $y_j'$ be the last point on $\gamma_j$ that is in $Z_L$.
Note that $d(x_j',x_j)/\operatorname{length}(\gamma_j)\to 0$ as $j\to\infty$,
because $L<L_1=\Lip_V u_1$,~\eqref{e.tight} holds and
$|u(x_j')-u(x_j)|\le L\,d(x_j',x_j)$.
 (If at some significant proportion of the length of $\gamma_j$ the function
$u$ does not change in speed very close to $L_1$, it will not have enough distance to catch up.)
Similarly, $d(y_j',y_j)/\operatorname{length}(\gamma_j)\to 0$.
Since $\sup_{x\notin Z_L}\Lip_x u_1\le L$ and $\sup_{x\notin Z_L}\Lip_x w\le L$
and $u_1=w$ on $Z_L$,
it follows that $|u_1(x_j)-w(x_j)|\le 2 \,L\,d(x_j,x_j')$
and similarly
$|u_1(y_j)-w(y_j)|\le 2 \,L\,d(y_j,y_j')$.
 Thus, $|w(x_j)-w(y_j)|/\operatorname{length}(\gamma_j)\to L_1$,
proving that $\Lip_{\partial V} w\ge \Lip_V w$, as needed.
\end{proof}

Our example is based on smoothing two different modifications of Aronsson's
$\infty$-harmonic function $G$ from~\eqref{e.Aro}.
For $L>0$ let $Z_L$ be the set of points in $\R^2\setminus\{0\}$
such that $|\nabla G|>L$.
Define the inner an outer radii of $Z_L$:
$ R_L:=\sup\{|z|:z\in Z_L\}$ and $r_L:=\inf \{|z|:0\ne z\notin Z_L\}$,
and note that $\lim_{L\to\infty}R_L=0$, while $r_L>0$ for every $L>0$.
Now fix some large $L$.
Let $u$ denote $G$ restricted to the annulus
$A_L:=\{z\in\R^2:r_L/2\le |z| \le 1\}$.
Let $u_j$, $j=1,2$, denote the AM function on $A_L$ whose boundary
values are equal to $u$ on the inner circle and are $j$ on the outer circle,
respectively.
Lemma~\ref{l.nocontin} implies that $u_1=u=u_2$ on $A_L\cap Z_L$
(provided that $L$ was chosen sufficiently large).
Shortly, we will show that there is a function $v$
defined on the unit disk that agrees
with $u$ in $A_L$ and such that $\Delta_\infty v=g$
(in the viscosity sense), where
$g$ is a Lipschitz function. Assuming this for the moment,
we can then define $v_j=u_j-j$ in $A_L$ and $v_j=v-j$
inside the disk $|z|<r_L/2$. Then $v_1$ and $v_2$
are both $0$ on the unit circle, and both
satisfy $\Delta_\infty v_j=g$, providing the required example.
It thus remains to prove

\begin{lemma}\label{l.smoothing}
For every $r_0>0$ there are Lipschitz functions $g,v:\R^2\to\R$
such that $v$ agrees with $G$ on $|z|>r_0$ and
$\Delta_\infty v=g$ holds in the viscosity sense.
\end{lemma}

\begin{proof}
Let
$$
a(\theta):=
\left[
\frac{\cos \theta\,(1-|\tan(\theta/2)|^{4/3})^2}{1 + |\tan(\theta/2)|^{4/3} +
  |\tan(\theta/2)|^{8/3}}\right]^{1/3}
$$
be the angle factor of $G$ in~\eqref{e.Aro}.
First, observe that $a(\theta+\pi)=-a(\theta)$ and $a(-\theta)=a(\theta)$.
Next, note that $a(\theta)>0$ when $\theta\in(-\pi/2,\pi/2)$.
We now verify that $a(\theta)$ is $C^\infty$ in a neighborhood of $\pi/2$.
Write
$$
a(\theta)=\cos \theta\,
\Bigl({1 + |\tan(\theta/2)|^{4/3} +
  |\tan(\theta/2)|^{8/3}}\Bigr)^{-1/3}
\Bigl(\frac{1-|\tan(\theta/2)|^{4/3}}{\cos\theta}\Bigr)^{2/3}.
$$
The first two factors are real-analytic near $\pi/2$.
Setting $\phi=|\tan(\theta/2)|^2/3$, which is real-analytic near $\theta=\pi/2$,
allows us to rewrite the last factor as
\begin{multline*}
\bigl(\cos{({\theta}/ 2)}\bigr)^{-4/3}\,
\Bigl(\frac{1-|\tan(\theta/2)|^{4/3}}{1-\tan^2(\theta/ 2)}\Bigr)^{2/3}=
\bigl(\cos{({\theta}/ 2)}\bigr)^{-4/3}\,
\Bigl(\frac{1-\phi^2}{1-\phi^3}\Bigr)^{2/3}
\\=
\bigl(\cos{({\theta}/ 2)}\bigr)^{-4/3}\,
\Bigl(\frac{1+\phi}{1+\phi+\phi^2}\Bigr)^{2/3},
\end{multline*}
which is also real-analytic near $\theta=\pi/2$.
We conclude that $a(\theta)$ is real-analytic
at every $\theta\notin \pi\,\Z$.

It is instructive to see that $G$ is $\infty$-harmonic in $\R^2\setminus\{0\}$.
Verifying $\Delta_\infty G=0$ away from the real line can be done by differentiation.
At points where $\theta=0, r>0$, we have $\nabla G/|\nabla G|=(-1,0)$.
At such points $\partial_x ^2 G=\partial_x^2 r^{-1/3}=4\,r^{-7/3}/9$, and indeed
$\Delta^+_\infty G=4\,r^{-7/3}/9$ there.
However,
\begin{equation}\label{e.aexpansion}
a(\theta)= 1- 16^{-1/3}\,\theta^{4/3}-\theta^2/6+O(\theta^3)\,,
\end{equation}
near $\theta=0$. This implies that there is no $C^\infty$ function $u$
that satisfies $u(r,0)=G(r,0)$ and $u\le G$ in a neighborhood of $(r,0)$.
Thus, $\Delta_\infty^-G=-\infty$.

We return to the proof of the lemma, and first
establish that the lemma holds when we only require that
$g$ be continuous, instead of Lipschitz.
In this case, the function $v$ may be written as
\begin{equation}\label{e.vdef}
v(r,\theta):=\phi(r)\,\bigl(\lambda(r)\,a(\theta)+(1-\lambda(r))\,\cos \theta\bigr)
\end{equation}
for appropriately chosen functions $\phi$ and $\lambda$.
The functions $\phi$ and $\lambda$ will be $C^\infty$, will satisfy
$\lambda(r)=\phi(r)=0$ in a neighborhood of $0$ and
$\phi(r)=r^{-1/3}$, $\lambda(r)=1$ for all $r\ge r_0/2$.
It follows that $v$ is Lipschitz and that $\Delta_\infty v$
is $C^\infty$ away from the $x$-axis and in a neighborhood of $0$.
It remains to check the behavior of $\Delta_\infty v$
at points on and close to the $x$-axis.
Based on~\eqref{e.aexpansion}, we may estimate $\Delta_\infty v$
for $\theta\ne 0$ close to $0$ and $r>0$:
\begin{equation}\label{e.dv}
\Delta_\infty v=
\phi''-\frac{4\,\lambda^3\,\phi^3}{81\,r^4\,(\phi')^2}
+\frac{P(r,\lambda,\phi,\phi',\lambda')}{r^6\,(\phi')^4}\,\theta^{2/3}+O(\theta)\,,
\end{equation}
where $P$ is some polynomial.
Let $r_1,r_2,r_3$ satisfy $0<r_3<r_2<r_1<r_0/2$.
Choose $\lambda(r)$ as a $C^\infty$ function that is $0$ for $r\le r_2$,
is $1$ for $r\ge r_1$, and is strictly monotone in $[r_2,r_1]$.
Choose $\phi(r)$ as a $C^\infty$ function that is
$0$ in $[0,r_3]$ and is  $r^{-1/3}$
for $r\ge r_2$.
When $r< r_2$, $\lambda=0$, and hence
$\Delta_\infty v$ is $C^\infty$ by~\eqref{e.vdef}.
In the range $r\in[r_2,r_1]$, we use~\eqref{e.dv}
 to conclude that
the limit as $\theta\to 0$ of $\Delta_\infty v$ exists.

We now argue that there is a continuous function $g$
such that
$\Delta_\infty v =g$ in the viscosity sense.
Suppose that $\eta$ is a $C^2$ function defined
in a small open  set $U$ containing the point $z_0=(r,0)$
such that the minimum of $\eta-v$ in $U$ occurs at $z_0$.
We will now perturb $\eta$.
For sufficiently small $s_1>0$ the function $\eta_1(z):=\eta+s_1\,|z-z_0|^2$
is a $C^2$ function defined in $U$ and $z_0$
is the unique minimum of $\eta_1-v$ in $U$.
If $s_2>0$ is chosen much smaller,
then the function $\eta_2(z)=\eta_1(z)+s_2\,y$
will still satisfy that the infimum of $\eta_2-v$
is attained in $U$. Moreover, that infimum
cannot be attained on the $x$-axis,
because $s_2\,y$ is zero on the $x$-axis,
$\nabla y$ is not parallel to the $x$-axis,
and $\nabla v$ exists.
Since $s_2$ and $s_1$ are arbitrarily small,
we conclude that $\Delta_\infty \eta(z_0)\ge \lim_z \Delta_\infty v(z)$
as $z$ tends to $z_0$ through point not on the $x$-axis.
Thus, $\Delta_\infty^+v(z_0) \ge \lim_z\Delta_\infty v(z)$.
Similarly, $\Delta_\infty^-v(z_0) \le \lim_z\Delta_\infty v(z)$.
This proves $\Delta_\infty v=g$, where $g$ is the continuous
extension of $\Delta_\infty v$ off the $x$-axis to the
whole plane.

If we want $g$ to be Lipschitz, we need to eliminate the
$\theta^{2/3}$ term in~\eqref{e.dv}.
Define
\begin{equation*}
\begin{aligned}
a_1(\theta) &:= -16^{1/3}\bigl( a(\theta)-\cos\theta - (1/12)\,\sin^2 2\theta\bigr)\,,\\
a_2(\theta) &:= \cos\theta + (1/8)\,\sin^2 2\theta\,,\\
a_3(\theta) &:= (1/4)\,\sin^2 2\theta\,.
\end{aligned}
\end{equation*}
Then near $\theta=0$ we have
\begin{equation*}
\begin{aligned}
a_1(\theta) &=  \theta^{4/3}+O(\theta^3)\,,
\\
a_2(\theta) &= 1+O(\theta^3)\,,
\\
a_3(\theta) &= \theta^2+O(\theta^3)\,.
\end{aligned}
\end{equation*}
Set
\begin{equation}\label{e.vdef2}
v(r,\theta):=\lambda_1(r) \, a_1(\theta)+\lambda_2(r)\,a_2(\theta)+\lambda_3(r)\,a_3(\theta)\,,
\end{equation}
where $\lambda_1,\lambda_2$ and $\lambda_3$ are to be chosen soon.
A calculation shows that near $\theta=0$ we have
\begin{multline*}
\Delta_\infty v =
\lambda_2''+\frac{64\,\lambda_1^3}{81\,r^4\,(\lambda_2')^2} +{}\\
{}- 16\,\lambda_1
\frac
{64\,\lambda_1^4 -162\,r^4\,\lambda_1'(\lambda_2')^3+27\,r^2\,\lambda_1\,(\lambda_2')^2\,\bigl(
3\,r^2\,\lambda_2''+3\,r\,\lambda_2'-10\,\lambda_3
\bigr)
}
{729\,r^6\,(\lambda_2')^4}
\,\theta^{2/3}+O(\theta)\,.
\end{multline*}
Of course, when $\lambda_1,\lambda_2$ and $\lambda_3$ are chosen so that
$v=G$, all the terms on the right hand side vanish.
Our goal is to choose these functions so that the
following holds:
\begin{enumerate}
\item  for $r\ge r_0$
we have $v=G$,
\item when $r$ is sufficiently small, we have $v=0$,
\item the $\theta^{2/3}$ term vanishes identically,
\item we don't have a blow up due to $\lambda_2'=0$,
and finally
\item the functions $\lambda_j$ are $C^\infty$,
say.
\end{enumerate}
This is not hard.
Fix $r_1,r_2,r_3,r_4$ such that $0<r_4<r_3<r_2<r_1<r_0/2$.
For $r>r_1$, we choose the $\lambda_j$ so that $v=G$.
In particular, $\lambda_2=r^{-1/3}$ in this range.
In the range $r\ge r_3$ we maintain $\lambda_2=r^{-1/3}$.
The function $\lambda_1$ is chosen
so that throughout $[r_3,r_2]$ we have $\lambda_1=\lambda_2'$,
while $\lambda_1$ does not vanish in $[r_3,r_1]$.
This is possible, because $\lambda_2'<0$ for $r\ge r_3$ and $\lambda_1=-16^{-1/3}<0$
at $r=r_1$. For every $r\in [r_3,r_1]$, $\lambda_3(r)$ takes the unique value
for which the $\theta^{2/3}$ term vanishes. Since $\lambda_1$ and $\lambda_2'$
do not vanish in the interval, it is clear that such a choice
for $\lambda_3$ is possible and $\lambda_3$ is $C^\infty$
provided that $\lambda_1$ and $\lambda_2$ are.
Throughout $r<r_2$ we take $\lambda_1=\lambda_2'$.
This allows us to simplify the expression for $\Delta_\infty v$ in this range, to obtain
$$
\Delta_\infty v
=
\lambda_2''+\frac {64 \,\lambda_2'}{81\,r^4}
+
16\,\frac{270\,r^2\,\lambda_3-(64+81\,r^3)\,\lambda_2'+ 81\,r^4\,\lambda_2''}{729\,r^6}\,\theta^{2/3}+O(\theta)\,.
$$
Now the denominators cannot vanish before $r=0$. Thus, $\lambda_3$ takes care to
make the $\theta^{2/3}$ term vanish while $\lambda_2$ evolves to become zero
throughout $r\le r_4$. This completes the proof.
\end{proof}

\section{Viscosity solutions and quadratic comparison in $\R^n$} \label{visc-comp}

In this section we prove Theorem~\ref{thm-visc-comp}, which states that in
bounded domains in Euclidean space $\R^n$,
$u$ is a viscosity solution of $\Delta_\infty u = g$
iff $u$ satisfies  $g$-quadratic comparison.
 The following lemma will be useful in that
proof. Let $D^2\phi(x)$ denote the Hessian matrix
$\bigl(\partial_i\partial_j\phi(x)\bigr)^{i=1,\ldots,n}_{j=1,\ldots,n}$.

\begin{lemma} \label{sandwich}
  Let $\phi$ be any real-valued function which is $C^2$ in a
  neighborhood of $x_0 \in \R^n$ and satisfies $\nabla \phi(x_0) \neq
  0$. Then for every $\delta >0$, there exist quadratic distance
  functions $\phi_1$ and $\phi_2$ such that:
  \begin{enumerate}
  \item $\phi_1(x_0) = \phi(x_0) = \phi_2(x_0)$;
     $\nabla\phi_1(x_0)=\nabla \phi(x_0)=\nabla\phi_2(x_0)$.
  \item $0<|\Delta_\infty \phi_i(x_0) - \Delta_\infty \phi(x_0)| <
    \delta$ for $i \in \{1,2\}$.  Also, as quadratic forms,
    $D^2\phi_2(x_0)$ strictly dominates $D^2\phi(x_0)$, which in turn
    strictly dominates $D^2\phi_1(x_0)$.
  \item Consequently, $\phi_1 < \phi < \phi_2$ at all points $x\ne x_0$
    in a neighborhood
    $\widetilde U$ of $x_0$.
  \item $\phi_1$ and $\phi_2$ are centered at $z_1$ and $z_2$, respectively, with
    $0<|z_i - x_0| < \delta$ for $i \in \{1,2\}$.
  \item In a neighborhood of $x_0$, the function $\phi_1$ is
      \decreasing in the distance from $z_1$ and $\phi_2$ is
      \increasing in the distance from $z_2$.
  \end{enumerate}
\end{lemma}

\begin{proof}
  We begin with preliminary calculations about the quadratic distance
  function $\phi_{a,b}(x) = a|x|^2 + b|x|$, centered at the origin.
  Fix some $x \neq 0$. Let $v = x/|x|$ be the unit vector in the
  $x$-direction and $v_O$ be any unit vector orthogonal to $v$.  Write
  $M = M(x) = D^2\phi_{a,b}(x)$.

  Then direct calculations show:

  \begin{enumerate}
  \item $\nabla \phi_{a,b}(x) = (2a |x| + b)\,v$.
  \item $v^T M v = 2a$.
  \item $v_O^T M v_O = 2a + b|x|^{-1}$.
  \item $v_O^T M v = 0$.
  \item $\Delta_\infty \phi_{a,b}(x) = 2a$ whenever $\nabla
    \phi_{a,b}(x) \neq 0$.
  \end{enumerate}
  Of course, it is enough to complete this calculation in dimension
  two when $x = (1,0)$ and $a=1$ and deduce the general case by symmetry.

  We now construct the $\phi_2$ described in the lemma (the
  construction of $\phi_1$ is similar).  We will take $z_2 = x_0 -
  \gamma \frac{\nabla\phi(x_0)}{|\nabla\phi(x_0)|}$ for some small
  value of $\gamma>0$ (specified below) and write $\phi_2(x) =
  \phi_{a,b}(x-z_2)+c$, where we first choose $a$ to be an arbitrary
  real with the property that $0 < 2 a - \Delta_\infty \phi(x_0) <
  \delta$, we then choose $b$ so that $\nabla \phi_2(x_0) = \nabla
  \phi(x_0)$ and we then choose $c$ so that $\phi_2(x_0) = \phi(x_0)$.
  We must now show that the $\phi_2$ thus constructed satisfies the
  requirements of the lemma.

  We can compute $b$ explicitly in terms of $a$, $x_0$,
  $\nabla\phi(x_0)$, $\gamma$, and $z_2$ using the relation
  $\nabla\phi_2(x_0) = (2a \gamma + b)\frac{x_0-z_2}{|x_0-z_2|} =
  \nabla \phi(x_0)$.  As $\gamma$ tends to zero, $b$ tends to $|\nabla
  \phi(x_0)|$.

  The description of $D^2 \phi_{a,b}(x)$ given above implies that if $v =
  (x_0-z_2)/|x_0-z_2|$ and $v_O$ is a unit vector orthogonal to $v$,
  and $M = D^2\phi_{a,b}(x_0)$, then
  \begin{enumerate}
  \item $v^T M v = 2a$.
  \item $v_O^T M v_O = 2a + b\gamma^{-1}$.
  \item $v_O^T M v = 0$.
  \end{enumerate}

  Note that $\chi:=v_O^T M v_O=2a+b\gamma^{-1}=\gamma^{-1}|\nabla \phi(x_0)|$
  tends to $\infty$ as $\gamma\searrow 0$.
  Let $M':=D^2\phi(x_0)$. We claim that by choosing $\gamma$
  sufficiently small we can make sure that $M$ dominates $M'$ as
  a quadratic form (or any other fixed quadratic form satisfying
  $v^TM'v<2a$, for that matter).
  Let $C:= \sup\{w_0^T M' w_1: w_0,w_1\in \R^n,\,|w_0|=|w_1|=1\}$
  and $a':=v^TM'v/2<a$.
If $w\in \R^n$ is nonzero, we may write
$w=\alpha\,v+w_O$, where $\alpha\in\R$ and $w_O$ is orthogonal to $v$.
Then
$$
w^TM'w
\le  2\,a'\,\alpha^2 + 2\,C\,\alpha\,|w_O|+ C\,|w_O|^2.
$$
On the other hand
$ w^TMw= 2\,a\,\alpha^2+\chi\,|w_O|^2$.
Since $2\,C\,\alpha\,|w_O|\le \frac{a-a'}2\,\alpha^2+C'\,|w_O|^2$
for some constant $C'=C'(C,a,a')$, the domination follows
from $\lim_{\gamma\searrow 0}\chi=\infty$.
This constructs $\phi_2$ with the required properties.
A similar construction applies to $\phi_1$.
\end{proof}

\begin{proof}[Proof of Theorem~\ref{thm-visc-comp}]
Let $x_0\in U$ and suppose that
$u$ satisfies $g$-quadratic comparison in a neighborhood of $x_0$.
Let $\phi$ be a $C^2$ real valued function defined in a neighborhood
of $x_0$. Suppose that $\nabla\phi(x_0)\ne 0$
and $\phi-u$ has a local minimum at $x_0$.
The above lemma implies that for every $\delta>0$ there is a quadratic distance
function $\phi_2(x)=a\,|x-z|^2+b\,|x-z|$ such that $\phi_2(x)-\phi(x)>\phi_2(x_0)-\phi(x_0)$
for all $x\ne x_0$ in some neighborhood of $x_0$,
$\Delta_\infty \phi_2<\Delta_\infty \phi+\delta$
and $\phi_2$ is \increasing in a neighborhood of $x_0$.
Since $\phi_2-u$ has a strict local minimum at $x_0$,
it follows by $g$-quadratic comparison that arbitrarily close
to $x_0$ there are points $x$ for which $a>g(x)/2$.
By continuity of $g$, this implies $a\ge g(x_0)/2$.
That is, $\Delta_\infty \phi_2(x_0)\ge g(x_0)$.
Since $\delta>0$ was arbitrary, this also implies
$\Delta_\infty\phi(x_0)\ge g(x_0)$.

Now we remove the assumption that $\nabla\phi(x_0)\ne 0$ and
assume instead that $\phi(x)=a\,|x-x_0|^2+o(|x-x_0|^2)$
as $x\to x_0$.
If $a\ge 0$, we may take $\phi_2(x)=(a+\delta)\,|x-x_0|^2$,
and the same argument as above gives
$\Delta_\infty \phi(x_0)\ge g(x_0)$.
Now suppose $a<0$ and let
$\delta\in (0,|a|)$. In this case, we have to modify the argument,
because $\phi_2$
is not \increasing in a neighborhood of $x_0$.
Recall that $\phi-u$ has a local minimum at $x_0$
and $\phi_2-\phi$ has a strict local minimum at $x_0$.
Therefore, $\phi_2-u$ has a strict local minimum at $x_0$.
Let $r>0$ be sufficiently small so that $\phi_2(x)-u(x)>\phi_2(x_0)-u(x_0)=-u(x_0)$
on $\{x\in\R^n: 0<|x-x_0|\le r\}$.
By continuity, for every $\eta\in \R^n$ with $|\eta|$
sufficiently small
$\phi_2(x+\eta)-u(x)>-u(x_0)$ for every $x\in \partial B_{x_0}(r)$.
Fix such an $\eta$ satisfying $0<|\eta|<r$, and let
$x_*$ be a point in $\overline{ B_{x_0}(r)}$
where $\phi_2(x+\eta)-u(x)$ attains its minimum.
Since
$$
\phi_2(x_0+\eta)-u(x_0)<-u(x_0)\le \inf\bigl \{ \phi_2(x+\eta)-u(x):x\in \partial B_{x_0}(r)\bigr\},
$$
it follows that $x_*\notin \partial B_{x_0}(r)$.
On the other hand, $x_*\ne x_0-\eta$,
because
$$
\phi_2\bigl((x_0-\eta)+\eta\bigr)-u(x_0-\eta)=- u(x_0-\eta)\ge -\phi_2(x_0-\eta)-u(x_0)>-u(x_0)\,.
$$
Therefore, $\nabla \phi_2(x_*)\ne 0$. The first case we have considered therefore
gives $g(x_*)\le \Delta_\infty \phi_2(x_*)=2\,(a+\delta)$.
Since $\delta\in (0,|a|)$ was arbitrary, continuity gives $g(x_0)\le 2\,a$, as required.
Thus, we conclude that $\Delta_\infty^+u\ge g$. By symmetry,
$\Delta_\infty^-u\le g$, and hence $u$ is a viscosity solution as claimed.

This completes the first half of the proposition.  For the converse,
suppose that $u$ is a viscosity solution of $\Delta_\infty u=g$.
Let $W\subset \overline W\subset U$ be open.
Suppose that $\phi(x)=a\,|x-z|^2+b\,|x-z|+c$ is a \increasing quadratic distance
function in $W$, $\phi >u$ on $\partial W$
and there is some $x_0\in W$ such that $u(x_0)>\phi(x_0)$.
For all sufficiently small $\delta>0$
we still have $\phi_\delta(x):=\phi(x)-\delta\,|x-z|^2>u$ on $\partial W$.
Let $x_*$ be a point in $\overline W$ where $\phi_\delta(x)-u(x)$ attains its
minimum.
Note that $x_*\notin\partial W$, since $\phi_\delta(x_0)-u(x_0)\le \phi(x_0)-u(x_0)<0$.
Consequently, $g(x_*)\le\Delta_\infty^+ u(x_*)\le \Delta_\infty\phi_\delta(x_*) = 2\,(a-\delta)$.
Thus, $2\,a>\inf_W g$. Therefore, $u$ satisfies $g$-quadratic comparison on $W$,
and the proof is complete.
\end{proof}

\section{Limiting trajectory} \label{limitingtrajectory}

\subsection{A general heuristic}
For every $\eps$, when both players play optimally in
$\eps$-tug-of-war, the sequence $\{x_k\}$ of points visited is random.
Do the laws of these random sequences, properly normalized, converge
in some sense to the law of a random continuous path as $\eps$ tends
to zero?

We give a complete answer in only a couple of simple cases.  However,
we can more generally compute the limiting trajectory when $u$ is
$C^2$ in a domain contained in $\R^n$ and the players move to
maximize/minimize $u$ instead of $u^\eps$; we suspect but cannot prove
generally that the limiting behavior will be the same when the players
use $u^\eps$.

Consider a point $x_0$ in the domain at which $\nabla u(x_0)\neq 0$.
If $\eps$ is small enough, then
$\nabla u\neq 0$ throughout the closed ball $\overline B_\eps(x_0)$, so the extrema
of $u$ on the closed ball $\overline B_\eps(x_0)$ lie on the surface
of the ball.  Then at any
such extremum $x$, by the Lagrange multipliers theorem, $x-x_0 = \lambda\nabla u(x)$
for some real $\lambda=O(\eps)$.
Since $u$ is $C^2$, we have $\nabla u(x) = \nabla u(x_0) + D^2 u(x_0)(x-x_0) + o(\eps)$,
where $D^2 u(x_0)=(\partial_i\partial_j
u(x_0))^{i=1,\ldots,n}_{j=1,\ldots,n}$.
We define
$\eta=\nabla u(x_0)$, $H = D^2 u(x_0)$, and $c=\eta^{T} H \eta/|\eta|^2$.
Then $x-x_0 = \lambda (\eta + H (x-x_0) + o(\eps))$, where $\lambda=O(\eps)$.
If we solve this with small $\lambda$, we find
$$x-x_0 = (I-\lambda H)^{-1}\lambda(\eta+o(\eps))  = \lambda \eta + \lambda^2 H \eta + o(\eps^2).$$
Then $\eps^2 = |x-x_0|^2 = |\eta|^2\lambda^2 + 2 c |\eta|^2 \lambda^3 + o(\eps^3)$,
and thus $\pm\eps = |\eta|\lambda + c|\eta|\lambda^2 + o(\eps^2)$, and
so $\lambda = \pm \eps/|\eta| - \eps^2 c/|\eta|^2 + o(\eps^2)$.  Hence,
$$x-x_0 = \pm \eps|\eta|^{-1} \eta + \eps^2|\eta|^{-2} (H- c I) \eta + o(\eps^2).$$
When $u$ is $C^2$ and has nonzero gradient in a domain, this suggests the SDE:
\begin{equation} \label{theSDE} dX_t =  r(X_t)\, dB_t + s(X_t)\, dt\,, \end{equation}
where $r(X_t) = |\nabla u(X_t)|^{-1} \nabla u(X_t)$ and $s(X_t)$ is equal to
$|\nabla u(X_t)|^{-2} D^2 u(X_t) \nabla u(X_t)$ minus its projection
onto $\nabla u(X_t)$ (so that $r(X_t)$ and $s(X_t)$ are always orthogonal).
Now It\^o's formula implies that, as expected,
$u(X_t) - \frac12\,\int_{s=0}^t \Delta_\infty u(X_s)\,ds$ is a martingale.  In particular, when
$u$ is infinity harmonic, $u(X_t)$ is a martingale.

All of the above analysis applies only when players make moves to
optimize $u$---instead of $u^{\eps}$, as they would do if they
were playing optimally.  This difference is what makes the calculation
of the limiting trajectory a heuristic, except in a few special cases
as described in the next subsection.

\subsection{Special cases}

The above analysis does apply to
optimally played games in a couple of simple cases for which $u =
u^{\eps}$.  Let $X = \R^2$ and let $Y\subset X$ be the complement
of a bounded set.
Then the following are infinity harmonic functions $u$ which satisfy
both $\Delta_{\infty} u=0$ and $\Delta_\infty^\eps u = 0$:

\begin{enumerate}
\item $u(v)$ is the distance from $v$ to a fixed convex set whose
$\eps$-neighborhood is contained in $Y$
\item $u(v)$ is the argument of $v$ on a $2\,\eps$-neighborhood of $X\setminus Y$,
 and defined arbitrarily elsewhere (here we assume that the $2\,\eps$-neighborhood of
 $X\setminus Y$ does not contain a simple closed curve surrounding $0$,
 so that the argument can be defined there).
\end{enumerate}
(Due to boundary errors the above need not hold when
$X \setminus Y$ is a bounded domain and $Y$ is its boundary.)
In the first case, $X_t$ is simply a Brownian motion along a
(straight) gradient flow line of $u$.  In the second case, $X_t$ is a
diffusion with drift in the $-X_t$ direction and diffusion of constant
magnitude orthogonal to $X_t$.

\medskip

Crandall and Evans \cite{MR1804769} have explored the following
question in some detail, and it has recently been answered
affirmatively by Savin \cite{savin} when $n=2$: Is every $\infty$-harmonic function
on a domain in $\R^n$ everywhere differentiable?  This question can be
rewritten as a question about the amount of variation of the optimal
direction of the first move (as a function of the starting point) in
$\eps$ tug-of-war.  An affirmative answer might be a step towards a
more complete analysis of the limiting game trajectories when $f=0$
and when $u$ is not smooth, since it would at least ensure that $r(X_t)$
is well-defined everywhere that the gradient is non-zero.

\section{Additional open problems}\label{sec:open}

\begin{enumerate}
\item
If $U\subset\R^n$ is open and bounded, $F:\partial U\to\R$ Lipschitz,
$g:\overline U\to \R$ Lipschitz,
is there a unique viscosity solution for
$\Delta_\infty u=g+c$ with boundary values given by $F$
for generic $c\in\R$. Here, generic could mean in the sense of
Baire category, or could mean almost every, or perhaps this is true
except for a countable set of $c$.
\item

Does Theorem \ref{deltainftyuisf} continue to hold if $F$ and $f$ are
merely continuous instead of uniformly continuous?  Can the $\inf |f|>0$
requirement be replaced with $f \geq 0$ (or $f \leq 0$) when $X$ has
finite diameter?  When solving
$\Delta_\infty u = g$ on bounded domains in $\R^n$, can the condition
that $g$ be continuous be replaced with a natural weaker condition (e.g.,
semicontinuity or piecewise continuity)?
\item
In Section~\ref{multiple-AM} we gave a triple $(X,Y,F)$ with $F$ positive
and Lipschitz for which the continuum value of tug-of-war is not the
unique AM extension of $F$.

Is there an example where $X$ is the closure of a connected open
subset of $\R^n$ and $Y$ is its boundary? In particular, let $X$ be
the set of points in $\R^2$ above the graph of the function
$|x|^{1+\delta}$, let $Y$ be the boundary of $X$ in $\R^2$, and let
$F(x,y)=y$ on $Y$. Is $y$ the only AM extension of $F$ to $X$? What
is the value of the corresponding game? (Note added in revision:
Changyou Wang and Yifeng Yu (personal communication) have recently
shown that $y$ is the only AM extension of $F$ to $X$.
The proof turned out to be rather simple.)

\item
Suppose that $u:U\to\R$, is Lipschitz and $g_1,g_2:U\to\R$ are continuous
on an open set $U\subset\R^n$ and that
$\Delta_\infty u=g_1$ as well as $\Delta_\infty u=g_2$, both in the
viscosity sense.
Does it follow that $g_1=g_2$?
(Note added in revision:
Yifeng Yu~\cite{YifengYu} proved this in the case $n=2$.)

\end{enumerate}

\section*{Acknowledgements}
We thank Alan Hammond and G\'abor Pete for useful discussions and
the referee for corrections to an earlier version.

\bibliographystyle{amsplain}
\bibliography{tug}

\end{document}